\documentclass[reqno,12pt,a4paper]{amsart}

\usepackage{amsmath,amsfonts,amsthm,amssymb,amsxtra}
\usepackage{tikz}
\usepackage{tikz-3dplot}
\usepackage{mathrsfs}
\usepackage{amssymb}
\usepackage{amsfonts}
\usepackage{latexsym}
\usepackage{amsthm}

\usepackage{graphicx}
\def\lb{\label}

\begin{document}

%%%%%%%%%% Some definitions %%%%%%%%%%

%%%%%%%% Equations, theorems %%%%%%%%%
\renewcommand{\theequation}{\arabic{section}.\arabic{equation}}
\theoremstyle{plain}
\newtheorem{theorem}{\bf Theorem}[section]
\newtheorem{lemma}[theorem]{\bf Lemma}
\newtheorem{corollary}[theorem]{\bf Corollary}
\newtheorem{proposition}[theorem]{\bf Proposition}
\newtheorem{definition}[theorem]{\bf Definition}
\newtheorem{remark}[theorem]{\it Remark}
%\theoremstyle{remark}
%\newtheorem{remark}[theorem]{\bf Remark}

%%%%% Alphabet %%%%%
\def\a{\alpha}  \def\cA{{\mathcal A}}     \def\bA{{\bf A}}  \def\mA{{\mathscr A}}
\def\b{\beta}   \def\cB{{\mathcal B}}     \def\bB{{\bf B}}  \def\mB{{\mathscr B}}
\def\g{\gamma}  \def\cC{{\mathcal C}}     \def\bC{{\bf C}}  \def\mC{{\mathscr C}}
\def\G{\Gamma}  \def\cD{{\mathcal D}}     \def\bD{{\bf D}}  \def\mD{{\mathscr D}}
\def\d{\delta}  \def\cE{{\mathcal E}}     \def\bE{{\bf E}}  \def\mE{{\mathscr E}}
\def\D{\Delta}  \def\cF{{\mathcal F}}     \def\bF{{\bf F}}  \def\mF{{\mathscr F}}
\def\c{\chi}    \def\cG{{\mathcal G}}     \def\bG{{\bf G}}  \def\mG{{\mathscr G}}
\def\z{\zeta}   \def\cH{{\mathcal H}}     \def\bH{{\bf H}}  \def\mH{{\mathscr H}}
\def\e{\eta}    \def\cI{{\mathcal I}}     \def\bI{{\bf I}}  \def\mI{{\mathscr I}}
\def\p{\psi}    \def\cJ{{\mathcal J}}     \def\bJ{{\bf J}}  \def\mJ{{\mathscr J}}
\def\vT{\Theta} \def\cK{{\mathcal K}}     \def\bK{{\bf K}}  \def\mK{{\mathscr K}}
\def\k{\kappa}  \def\cL{{\mathcal L}}     \def\bL{{\bf L}}  \def\mL{{\mathscr L}}
\def\l{\lambda} \def\cM{{\mathcal M}}     \def\bM{{\bf M}}  \def\mM{{\mathscr M}}
\def\L{\Lambda} \def\cN{{\mathcal N}}     \def\bN{{\bf N}}  \def\mN{{\mathscr N}}
\def\m{\mu}     \def\cO{{\mathcal O}}     \def\bO{{\bf O}}  \def\mO{{\mathscr O}}
\def\n{\nu}     \def\cP{{\mathcal P}}     \def\bP{{\bf P}}  \def\mP{{\mathscr P}}
\def\r{\rho}    \def\cQ{{\mathcal Q}}     \def\bQ{{\bf Q}}  \def\mQ{{\mathscr Q}}
\def\s{\sigma}  \def\cR{{\mathcal R}}     \def\bR{{\bf R}}  \def\mR{{\mathscr R}}
                \def\cS{{\mathcal S}}     \def\bS{{\bf S}}  \def\mS{{\mathscr S}}
\def\t{\tau}    \def\cT{{\mathcal T}}     \def\bT{{\bf T}}  \def\mT{{\mathscr T}}
\def\f{\phi}    \def\cU{{\mathcal U}}     \def\bU{{\bf U}}  \def\mU{{\mathscr U}}
\def\F{\Phi}    \def\cV{{\mathcal V}}     \def\bV{{\bf V}}  \def\mV{{\mathscr V}}
\def\P{\Psi}    \def\cW{{\mathcal W}}     \def\bW{{\bf W}}  \def\mW{{\mathscr W}}
\def\o{\omega}  \def\cX{{\mathcal X}}     \def\bX{{\bf X}}  \def\mX{{\mathscr X}}
\def\x{\xi}     \def\cY{{\mathcal Y}}     \def\bY{{\bf Y}}  \def\mY{{\mathscr Y}}
\def\X{\Xi}     \def\cZ{{\mathcal Z}}     \def\bZ{{\bf Z}}  \def\mZ{{\mathscr Z}}
\def\O{\Omega}

\newcommand{\gA}{\mathfrak{A}}
\newcommand{\gB}{\mathfrak{B}}
\newcommand{\gC}{\mathfrak{C}}
\newcommand{\gD}{\mathfrak{D}}
\newcommand{\gE}{\mathfrak{E}}
\newcommand{\gF}{\mathfrak{F}}
\newcommand{\gG}{\mathfrak{G}}
\newcommand{\gH}{\mathfrak{H}}
\newcommand{\gI}{\mathfrak{I}}
\newcommand{\gJ}{\mathfrak{J}}
\newcommand{\gK}{\mathfrak{K}}
\newcommand{\gL}{\mathfrak{L}}
\newcommand{\gM}{\mathfrak{M}}
\newcommand{\gN}{\mathfrak{N}}
\newcommand{\gO}{\mathfrak{O}}
\newcommand{\gP}{\mathfrak{P}}
\newcommand{\gQ}{\mathfrak{Q}}
\newcommand{\gR}{\mathfrak{R}}
\newcommand{\gS}{\mathfrak{S}}
\newcommand{\gT}{\mathfrak{T}}
\newcommand{\gU}{\mathfrak{U}}
\newcommand{\gV}{\mathfrak{V}}
\newcommand{\gW}{\mathfrak{W}}
\newcommand{\gX}{\mathfrak{X}}
\newcommand{\gY}{\mathfrak{Y}}
\newcommand{\gZ}{\mathfrak{Z}}

\def\ve{\varepsilon} \def\vt{\vartheta} \def\vp{\varphi}  \def\vk{\varkappa}

\def\Z{{\mathbb Z}} \def\R{{\mathbb R}} \def\C{{\mathbb C}}  \def\K{{\mathbb K}}
\def\T{{\mathbb T}} \def\N{{\mathbb N}} \def\dD{{\mathbb D}} \def\S{{\mathbb S}}
\def\B{{\mathbb B}}

%%%%% Arrows %%%%%

\def\la{\leftarrow}              \def\ra{\rightarrow}     \def\Ra{\Rightarrow}
\def\ua{\uparrow}                \def\da{\downarrow}
\def\lra{\leftrightarrow}        \def\Lra{\Leftrightarrow}
\newcommand{\abs}[1]{\lvert#1\rvert}
\newcommand{\br}[1]{\left(#1\right)}

\def\lan{\langle} \def\ran{\rangle}

%%%%% Typography %%%%%

\def\lt{\biggl}                  \def\rt{\biggr}
\def\ol{\overline}               \def\wt{\widetilde}
\def\no{\noindent}

%%%%% Math signs %%%%%

\let\ge\geqslant                 \let\le\leqslant
\def\lan{\langle}                \def\ran{\rangle}
\def\/{\over}                    \def\iy{\infty}
\def\sm{\setminus}               \def\es{\emptyset}
\def\ss{\subset}                 \def\ts{\times}
\def\pa{\partial}                \def\os{\oplus}
\def\om{\ominus}                 \def\ev{\equiv}
\def\iint{\int\!\!\!\int}        \def\iintt{\mathop{\int\!\!\int\!\!\dots\!\!\int}\limits}
\def\el2{\ell^{\,2}}             \def\1{1\!\!1}
\def\sh{\sharp}
\def\wh{\widehat}
\def\bs{\backslash}
\def\na{\nabla}
%%%%% Math operations %%%%%

\def\sh{\mathop{\mathrm{sh}}\nolimits}
\def\all{\mathop{\mathrm{all}}\nolimits}
\def\Area{\mathop{\mathrm{Area}}\nolimits}
\def\arg{\mathop{\mathrm{arg}}\nolimits}
\def\const{\mathop{\mathrm{const}}\nolimits}
\def\det{\mathop{\mathrm{det}}\nolimits}
\def\diag{\mathop{\mathrm{diag}}\nolimits}
\def\diam{\mathop{\mathrm{diam}}\nolimits}
\def\dim{\mathop{\mathrm{dim}}\nolimits}
\def\dist{\mathop{\mathrm{dist}}\nolimits}
\def\Im{\mathop{\mathrm{Im}}\nolimits}
\def\Iso{\mathop{\mathrm{Iso}}\nolimits}
\def\Ker{\mathop{\mathrm{Ker}}\nolimits}
\def\Lip{\mathop{\mathrm{Lip}}\nolimits}
\def\rank{\mathop{\mathrm{rank}}\limits}
\def\Ran{\mathop{\mathrm{Ran}}\nolimits}
\def\Re{\mathop{\mathrm{Re}}\nolimits}
\def\Res{\mathop{\mathrm{Res}}\nolimits}
\def\res{\mathop{\mathrm{res}}\limits}
\def\sign{\mathop{\mathrm{sign}}\nolimits}
\def\span{\mathop{\mathrm{span}}\nolimits}
\def\supp{\mathop{\mathrm{supp}}\nolimits}
\def\Tr{\mathop{\mathrm{Tr}}\nolimits}
\def\BBox{\hspace{1mm}\vrule height6pt width5.5pt depth0pt \hspace{6pt}}
\def\where{\mathop{\mathrm{where}}\nolimits}
\def\as{\mathop{\mathrm{as}}\nolimits}

%%%%%%%%%%%%% specialities %%%%%%%%%%%%%%

\newcommand\nh[2]{\widehat{#1}\vphantom{#1}^{(#2)}}
%{{\mathop{#1}\limits^\wedge}\vphantom{#1}^{(#2)}}
\def\dia{\diamond}

\def\Oplus{\bigoplus\nolimits}

%%%%%%%%%%% End of definitions %%%%%%%%%%

\def\qqq{\qquad}
\def\qq{\quad}
\let\ge\geqslant
\let\le\leqslant
\let\geq\geqslant
\let\leq\leqslant
\newcommand{\ca}{\begin{cases}}
\newcommand{\ac}{\end{cases}}
\newcommand{\ma}{\begin{pmatrix}}
\newcommand{\am}{\end{pmatrix}}
\renewcommand{\[}{\begin{equation}}
\renewcommand{\]}{\end{equation}}
\def\eq{\begin{equation}}
\def\qe{\end{equation}}
\def\[{\begin{equation}}
\def\bu{\bullet}

\title[{ Eigenvalue  bounds  for Stark operators with complex potentials}]
{ Eigenvalue  bounds  for Stark operators  with complex potentials}

\date{\today}

\author[Evgeny Korotyaev]{Evgeny  Korotyaev}
\address{E. Korotyaev, Saint-Petersburg State University, Universitetskaya nab.
7/9, St. Petersburg, 199034, Russia; \ korotyaev@gmail.com, \
e.korotyaev@spbu.ru}

\author[Oleg Safronov]{Oleg Safronov}
\address{O. Safronov, Department of Mathematics and Statistics, UNCC, Charlotte, NC, USA;\  osafrono@uncc.edu}

\subjclass[2010]{35P15   (47A75)} \keywords{Stark operators, complex
potentials, eigenvalue  estimates }

\thanks{{\it Acknowledgments.}  E. Korotyaev was supported
by the RSF grant  No. 18-11-00032 \\ ${\,\,\,\,\,}$  The authors are grateful to Natalia Saburova  for her help with the manuscript}

\begin{abstract}
\no We consider the 3-dimensional Stark operator  perturbed  by a  complex-valued  potential. We
obtain an  estimate  for the  number of eigenvalues of this operator as well as  for   the  sum  of imaginary parts  of eigenvalues   situated  in the upper  half-plane.
\end{abstract}

\maketitle

%\vskip -3.25cm
\section {Introduction and main results}
\setcounter{equation}{0}

 Let $H_0$ be   the  free Stark
operator
\[
\label{a.1} H_0=-{\Delta}+ x_1,
\]  acting in  the space
$L^2(\R^3)$.
In the formula  above,  $x_1$ denotes the function  whose  value  at a point $x\in \R^3$ coincides  with  the first  coordinate of  $x$.  Since $H_0$  is an unbounded  operator, one  has to specify   its  domain   of definition. For this purpose we  simply mention that  $H_0$ is essentially selfadjoint on $C_0^\infty({\mathbb R}^3)$.
We study the spectral properties  of  the operator $$H=H_0+V,$$
where the potential $V$ is  a  bounded  complex-valued      function, satisfying
\[
\lb{dV} %\int_{\R^3}|V|^2dx<\infty, \qquad
\int_{\R^3}|V(x)|^{r}dx<\iy,\qquad \text{for\,\,some}\quad  r>0.
\]
While   the interest  of   mathematicians  in the theory of  non-selfadjoint operators of this type   is quite   new,  Stark operators  with real potentials  have been studied thoroughly in mathematical physics   for   a  long time.
Among the classical results applicable to   the self-adjoint case  are    the theorems   of   Avron and Herbst \cite{AH77}   who considered
 scattering for the pair  of  operators $H$ and $H_0$ in the case where $V$ is  a short-range
potential.   In particular, it  was established that the spectrum of  $H$ is purely absolutely
continuous and covers the real line $\R$ (besides
\cite{AH77}, see  Herbst \cite{H77}).    It  was  proved in  \cite{AH77}, \cite{H77} and  \cite{Y81}  that   for a  short-range potential $V$, the wave
operators
$$
\Omega_\pm=s-\lim e^{itH}e^{-itH_0}, \qqq \text{as} \qqq t\to \pm\iy,
$$
exist and  are unitary.  Further development  of the  methods   used   to  study  Stark operators  led to  the  theory of
scattering of  several particles   in an external constant
electric field (see  the papers  \cite{HMS96}  and \cite{K87}).
Besides the  results   related  to the  scattering   theory, the   mathematical literature  on   Stark operators contains   numerous statements about the   distribution of   resonances  in the  models involving  a  constant  electric field.
 Here, we  only mention  the article
\cite{H79} and the recent paper \cite{K17}, which can be  also  used  for  finding  other   relevant  references.

Let us now  describe  the  main results  of the  present paper  devoted to  the non-selfadjoint case.
 Under  the condition  \eqref{dV},  the spectrum $\sigma(H)$ of $H$  coincides  with the union of  the line  ${\mathbb R}$ and the discrete set of complex eigenvalues  that   might accumulate  only to real  points.
 We denote by $\{\l_j\}_{j=1}^{\iy}$ the sequence of   eigenvalues  of $H$
in $\C\setminus \R$   enumerated  in an  arbitrary order.
   The  number of times  an  eigenvalue  appears in the sequence $\{\l_j\}_{j=1}^{\iy}$ coincides   with the  algebraic  multiplicity of  the eigenvalue.
We  will   show that  the  condition \eqref{dV} with $r<2$  guarantees  that
$$
\sum_j |{\rm Im}\,  \lambda_j|<\infty.
$$

\begin{theorem}\label{main} Let $4<p\leq5$.  Let $V$ be a  bounded complex-valued  function  satisfying    the condition \eqref{dV} with $r=p/(p-2)$.
Then the eigenvalues $\l_j$ of the operator $H$  obey the estimate
\[\lb{LTg=1}
\sum_j |\Im \l_j|\leq
C_p
\Bigl[\Bigl(\int_{\R^3}|V|^{p/2}dx\Bigr)^{2}+\Bigl(\int_{\R^3}|V|^{p/(p-2)}dx\Bigr)^{p-2}\Bigr].
\]
%%%%%%%%%%%%
%%%%%%%%%%%+\Bigl(\int_{\R^3}(1+|x_1|)^{3/4}|V|^{3/2}dx\Bigr)^{2/3}
The  constant $C_p>0$  in this inequality is independent of $V$.
\end{theorem}

\medskip

As  a  consequence of  the method used  in the proof of Theorem~\ref{main} we  will obtain the following  statement,  where  $V$  might decay slower than a potential satisfying \eqref{dV}  with $r<2$.
\begin{theorem}\lb{p>1} Let $q>1$  and $4<p<q+3.$
Let $V$ be a bounded complex-valued  potential such that $$\int_{{\mathbb R}^3}|V(x)|^{p/2}dx<\infty.$$
Then  the spectrum of  $H$ is discrete  in $\C\setminus \R$ and  the  eigenvalues $\lambda_j$   of the operator $H$   satisfy  the  estimate
\[
\lb{conseq8}
 \sum_j
 |{\rm Im}\, \lambda_j|^{q}\leq
C_{p,q} \Bigl(\int_{{\mathbb R}^3}|V(x)|^{p/2}dx\Bigr)^{2q/(p-3)}
\]
 with   a   positive constant $C_{p,q}>0$  depending only  on $p$ and $q$.
\end{theorem}

 \bigskip

{\bf Remark.} Estimates  \eqref{LTg=1}  and \eqref{conseq8} also  hold  for  eigenvalues  of the  operator $-\D+V$ in $d=3$.

\bigskip

The  next theorem gives a  very interesting  bound  on the number of eigenvalues of $H$  in the  half-plane $\{\l \in {\mathbb C}:\,\,{\rm Re}\, \l<\a\}$  under  the condition
\[\lb{1.3}
\int_{\R^3}   |V(x)|^{p/2}  \r(x)dx<\infty\]  where the  weight  $\r$ is the exponentially growing  as $x_1\to -\iy$ function,   given by
\[\lb{P=}
\r(x)=(1+e^{-p x_1/2})(1+|x_1|)^2\quad \text{and}\quad p>5.\]
Functions $V$ satisfying   the condition \eqref{1.3}  decay  exponentially   fast  in some integral sense   in the direction of the  negative $x_1$-axis. However, such potentials might  decay  slowly  in other  directions.   For instance,
any function   obeying
$$
|V(x)|\leq \frac{C}{(1+e^{-x_1})(1+|x|)^s},\qquad s>0,\,\,C>0,
$$
satisfies \eqref{1.3} with $p>10/s$. While the usual Schr\"odinger operator $-\D+V$ perturbed  by such a potential might have infinitely
many non-real eigenvalues  in $\{\l \in {\mathbb C}:\,\,{\rm Re}\, \l<\a\}$,  our  theorem says that the number of   non-real eigenvalues  of  the  Stark operator $H$ in this  half-plane is  still finite.

\begin{theorem}\label{main2} Let  $\d>0$  and $p>5$.  Let $V$ be a bounded  complex-valued  function satisfying \eqref{1.3}  with $\r(x)$  defined   by \eqref{P=}  and let $\a>0$.
Then the number $N(\alpha)$  of non-real  eigenvalues  of the operator $H$ situated in the  half-plane $\{\l \in {\mathbb C}:\,\,{\rm Re}\, \l<\a\}$  obeys the estimate
\[\lb{N<V}
N(\a)\leq  C_{\a,p,\d}
\Bigl(\int_{\R^3}  |V(x)|^{p/2}\r(x) dx\Bigr)^{2(1+\d)},
\]
where
$$ C_{\a,p,\d}= C_{p,\d} \min_{\varepsilon>0}  \Bigl[ \varepsilon^{-2} e^{(1+\d)p(\a+\varepsilon)}\Bigl( \frac{\alpha+\varepsilon}{\varepsilon^{1+2\d}}+ (1+\ve^2)e^{2(1+\d)p\varepsilon^2} \Bigr)\Bigr],$$  and
$C_{p,\d}>0$ is a  constant  that depends  only on $p$ and $\d$.

\end{theorem}

%%%%%%%%%%%%%%%%%%%%%%%%%%%%%%%%%%%%%%%%%%%%%%%%%%%%%%%%%%%%

\setlength{\unitlength}{1.0mm}
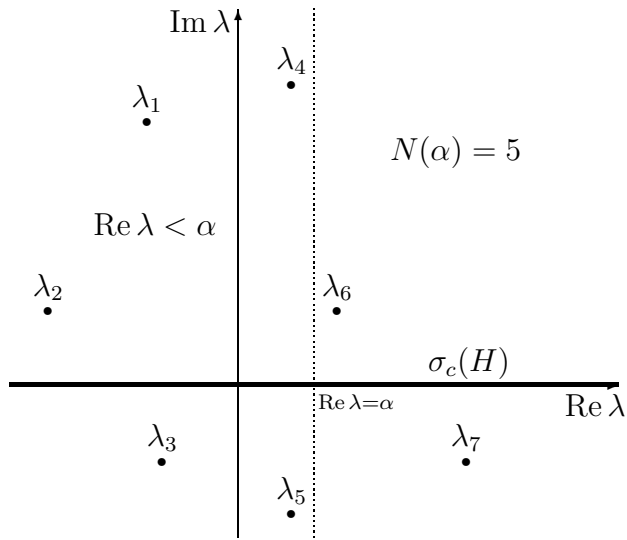
\begin{figure}[h]
\centering
\unitlength 1.0mm % = 2.845pt
%\linethickness{0.4pt}
%\ifx\plotpoint\undefined\newsavebox{\plotpoint}\fi % GNUPLOT compatibility
\begin{picture}(80,80)
\put(0,20){\vector(1,0){80.00}} \put(30,0){\vector(0,1){70.00}}

\qbezier[100](40,0)(40,35)(40,70)
\put(40.5,17){$\scriptstyle\Re\l=\a$} \put(21.0,67){$\Im\l$}
\put(73,16){$\Re\l$} \put(11.0,40){$\Re\l<\a$}
\put(50.0,50){$N(\a)=5$} \put(55.0,22){$\s_c(H)$}

\put(37,60){\circle*{1}} \put(35,62){$\l_4$}
\put(18,55){\circle*{1}} \put(16,57){$\l_1$} \put(5,30){\circle*{1}}
\put(3,32){$\l_2$} \put(43,30){\circle*{1}} \put(41,32){$\l_6$}
\put(20,10){\circle*{1}} \put(18,12){$\l_3$}
\put(60,10){\circle*{1}} \put(58,12){$\l_7$} \put(37,3){\circle*{1}}
\put(35,5){$\l_5$}

\linethickness{0.8pt} \put(0,20.3){\line(1,0){80.00}}
\end{picture}

%\vspace{-10mm}
\caption{\footnotesize   Eigenvalues of $H$}
%\label{ffS'}
\end{figure}

\medskip

\medskip

{\bf Remarks.} A  similar   estimate holds  for the number of resonances  of  $H$ contained  in a    region $\{\l\in \C:\,\, -\iy<{\rm Re}\, \l<\a,\,\,\,-\b<{\rm Im}\,\l\leq 0\}$   (the constant in such an estimate   depends  on the region).

\medskip

One can combine this theorem with  the fact  that all eigenvalues  are  situated  in a disk  of a finite radius, to obtain an estimate for the total number of non-real  eigenvalues.
\begin{theorem} \lb{individual}1)  Let $V\in L^\infty({\mathbb R}^3)$.
 There  exists  a universal    constant  $C>0$, such that
all  eigenvalues $\lambda\in {\mathbb C}\setminus \R$ of the operator $H$ are  situated in the  disk
\[\lb{ball}
|\l|\leq C\Bigl(\int_{{\mathbb R}^3}(1+|x|)^4|V(x)|\,dx+\Bigl(\int_{\R^3}|V|^2dx\Bigr)^{1/2}\Bigr)^4.
\]
In particular, the conditions
\eqref{1.3} and \eqref{ball}
 imply  that the total number $N$  of non-real eigenvalues of the operator $H$ is  finite
$$
\eqref{1.3}\,\, \text{and} \,\,\eqref{ball} \implies  N<\infty
$$ and coincides   with
$N(\alpha)
$,
where  $\a$  equals  the right  hand side of  \eqref{ball}.

2) If $H$  is the  Stark operator perturbed  by a potential  \[\lb{1.9}V\in L^{q/2}(\R^3)\cap  L^\infty(\R^3),\text{where}\quad q<3,\] then  non-real  eigenvalues $\l_j$ of $H$ are contained in   a disk of a finite radius.
In particular, if $V$ satisfies both  hypothesis \eqref{1.3} and \eqref{1.9},  then $H$ has  finitely  many eigenvalues in $\C\setminus \R$
$$
\eqref{1.3}\,\, \text{and} \,\,\eqref{1.9} \implies  N<\infty.
$$
\end{theorem}

\bigskip

{\bf Remark.}  In the same way, one can show that the condition $V\in L^{d/2}( \R^d)$ with $d\geq3$  implies  that all  non-real eigenvalues $z_j$ of $-\Delta+V$
are  contained in a disk of a finite radius  and  $\sum_j|{\rm Im }\sqrt{z_j}|<\infty$  (see \cite{FrSab}).
\bigskip

Let us say a  couple of  words  about
 our  approach to the problem.   It is   well known that  eigenvalues  of   most   important differential   operators   can be  described  as  zeros  of   the corresponding perturbation determinants,
which depend  analytically  on the spectral parameter.  The  latter observation allows  one   to  turn the analysis of eigenvalues into the study of zeros of analytic  functions.
Similar  ideas were  successfully used in  the paper    \cite{FrSab} by Frank and Sabin  for the study of the eigenvalues of the  Schr\"odinger operator perturbed  by a decaying  potential.
Among other related papers are  the articles \cite{Frank}, \cite{FLS}.  The problem pertaining to  the Stark operator is however more  complicated  compared  to the   one  involving the usual  Schr\"odinger equation,
simply  because  the  free Stark operator is  not diagonalized  by the  Fourier transformation.

We  would also  like to point out that the approach    based on the   study of the perturbation determinant  is not the only  method    of  obtaining   eigenvalue estimates in the  non-selfadjoint case.
For instance,    the authors of \cite{FLLS}  and \cite{LapSaf}   use a  completely  different technique to  estimate eigenvalues of a Schr\"odinger operator.

\bigskip

Our present paper has a complicated structure.
The proof of Theorem~\ref{main} is given in Sections 3 and 4.
Theorem~\ref{p>1} will be proved in Sections 5, 6 and 7.  The  two following  Sections 8 and  9  contain a proof of  Theorem~\ref{main2}  establishing a  bound  for the  number of non-real eigenvalues
contained in the half-plane $\{\l\in\C:\,\,{\rm Re\,\l}<\a\}$.
The   estimate \eqref{ball} of the radius of the  disk  containing all eigenvalues  will be   justified in the  last  Section 10.

\bigskip

 {\bf Notations.} We denote by $C$
various possibly different constants whose values are irrelevant.  The upper half-plane $\{\l\in\C: {\rm Im} \l>0\}$ will be  denoted by the symbol $\C_+$.
 By $\cB$ and ${\mathfrak S}_\iy$  we denote the classes of
bounded and compact operators, respectively. The symbols  $\mathfrak S_1$ and $\mathfrak S_2$  are used   to denote  the trace class and the Hilbert-Schmidt class  equipped with the norms
$\|\cdot \|_{\mathfrak S_1}$ and $ \|\cdot \|_{\mathfrak S_2}$, respectively.  More  generally,
$\mathfrak S_p$ denotes  the  class  of compact  operators $K$   obeying
$$
||K||^p_{\mathfrak S_p}={\rm tr}\Bigl(K^*K\Bigr)^{p/2}<\infty,\qquad p\geq1.
$$
Note that  if $K\in {\mathfrak S}_p$  for  some $p\geq1$, then $K\in {\mathfrak S}_q$  for $q>p$  and
$$
||K||_{{\mathfrak S}_q}\leq ||K||_{{\mathfrak S}_p}.
$$
For a  self-adjoint operator $T=T^*$ the symbol $E_T(\cdot)$ denotes its  (operator-valued) spectral measure.

\bigskip

\section{Preliminaries}

Very often,  eigenvalues  of closed   operators     can be described as zeros  of   analytic  functions.
The latter circumstance   allows  one  to use   known results on the distribution  of zeros of holomorphic   functions   to obtain   bounds on the eigenvalues of a given operator. In particular,
the
eigenvalues  of $H$ coincide with    zeros  of  the  so called   perturbation determinant $D_n(\l)$, which depends  analytically on $\l$.

The definition of $D_n(\l)$ requires that we find  two functions  $W_1$  and $W_2$  having the properties
$$
V=W_2W_1,\qquad |W_1|=|W_2|,
$$
and  set
\[
Y_0(\l)=W_1\, R_0(\l)\,
W_2, \qqq  R_0(\l)=(H_0-\l)^{-1},\qqq  \l \in {\C}\setminus \R.
\]
 The condition \eqref{dV} implies that
 $Y_0(\l)$ is an ${\frak S}_{2r}$- operator  whenever $r>3/2$  and   $ {\rm Im} \l \ne 0$.
Therefore, we can
define the   determinants
\[
\label{a.2} D_n(\l)=\det_n (I +Y_0(\l)),\qqq \l\in \C\setminus \R,
\]
for integer $n\geq 2r$.  The standard  way to  describe
$\det_n(I+K)$  in terms  of   eigenvalues $z_j$ of a  compact  operator $K\in {\mathfrak S}_n$ is   to  define  it as
$$
\det_n(I+K)=\prod_j(1+z_j)\exp\Bigl(\sum_{m=1}^{n-1}\frac{(-1)^mz_j^m}{m}\Bigr),\qquad n\geq 2;
$$
$$
\det(I+K)=\prod_j(1+z_j),\qquad  n=1.
$$

The following  relations   can be  found in Section 3 of the book \cite{S05}.
  If  $X, Y\in \cB$ and both products $XY$,
$YX$ belong  to $\mathfrak S_n$, then
\[
\lb{2.2} \det_n (I+ XY)=\det_n (I+YX).\]
 The\ mapping   $X\to \det (I+ X) $  is\  continuous on  $\mathfrak S_1$, which is  guaranteed  by the inequality
\[
\lb{B12}
|\det (I+ X)-\det (I+Y)|\le \|X-Y\|_{\mathfrak S_1}e^{1+\|X\|_{\mathfrak S_1}+\|Y\|_{\mathfrak S_1}}.
\]
Moreover,  there  exists a  constant $C_n>0$ (see \cite{Ha}) depending  only on $n$  such that
\[
\lb{B1}
|\det_n (I+ X)|\le e^{C_n\|X\|_{{\mathfrak  S}_n}^n},\qquad \forall X\in {\mathfrak S}_n.
\]
While the    inequality in the proposition below is  less  known compared to \eqref{B1},  it is  still  a very  useful  estimate  of the  $n$-th determinant.
\begin{proposition}  Let $n\geq2$. Then
 for any $n-1\leq p\leq n$,  there  exists a  constant $C_{p,n} >0$ depending  only on $p$ and $n$  such that
\[
\lb{det<p}
|\det_n (I+ X)|\le e^{C_{p,n}\|X\|_{{\mathfrak  S}_p}^p},\qquad \forall X\in {\mathfrak S}_p,\qquad n\geq2.
\]
\end{proposition}

{\it Proof.} We need to show that
$$
\ln(|\det_n (I+ X)|)\leq C_{p,n} \|X\|_{{\mathfrak  S}_p}^p
$$
with some constant $C_{p,n}>0.$ For that purpose, it is  sufficient  to prove that
\[\lb{dokazat2.17}
{\rm Re}\,\ln (1+z)+{\rm Re}\,\sum_{m=1}^{n-1}\frac{(-1)^mz^m}{m}\leq C_{p,n} |z|^p
\]
for all $z\in \C.$  The  inequality  \eqref{dokazat2.17}  is obvious  for very large and  very small  $|z|$.
Therefore,  it holds  for  all $z$  lying  outside of a small neighborhood of the point $-1$.
On the other  hand, the left hand  side  of \eqref{dokazat2.17}  is negative if $z$ is  sufficiently close to $-1$. Consequently, it holds  everywhere. $\BBox$

\bigskip

If  an  operator-valued function $X :\O\to \mathfrak S_1$ is  analytic on a domain $\O\ss\C$ and $(I+X(z))^{-1}\in \cB$ for all $z\in
\O$, then the function $F(z)=\det (I+X (z))$ is  also analytic and its  derivative satisfies the  relation
\[
\lb{2.3} F'(z)= F(z)\Tr\Bigl( (I+X (z))^{-1}X '(z)\Bigr),\ \   \ \ \ z\in \O.
\]
Similarly, if  an  operator-valued function $X :\O\to \mathfrak S_n$, (here,  $n\geq2$) is  analytic on a domain $\O\ss\C$ and $(I+X(z))^{-1}\in \cB$ for all $z\in
\O$, then the function $F(z)=\det_n (I+X (z))$ is analytic and its  derivative equals
\[
\lb{2.3b} F'(z)= F(z)\Tr\Bigl(\bigl( (I+X (z))^{-1}-\sum_{j=0}^{n-2} (-1)^j  X^j\bigr)X '(z)\Bigr),\ \   \ \ \ z\in \O.
\]

\bigskip

Let $n\geq2r$ be  integer. We will show that  if $V$  is  a bounded function satisfying  \eqref{dV} with $r>3/2$, then  $Y_0(\l)$  is an  ${\frak S}_{2r}$-operator  for
all $\l\in{\mathbb C}\setminus {\mathbb R}$. The latter condition  implies  that
the function $D_n(\l)$ is
analytic on the  open  domain $  {\mathbb C}\setminus {\mathbb R}$.
The following  statement is  known as  the Birman-Schwinger principle (for more detailed  description,  see  \cite{FLS}).

\begin{lemma} Let $V\in L^\infty(\R^3)$  satisfy \eqref{dV}  with $r>3/2$.  Let $n\geq2r$ be  integer.
The point $\l\in \C\setminus \R$ is an eigenvalue  of the operator   $H$  if and  only  if  $\l$ is a  zero of  $D_n(\l)$.
The  algebraic  multiplicity of each  eigenvalue $\l\in \sigma (H)\setminus \R$  coincides with the multiplicity of   the corresponding  zero of  the  function  $D_n(\cdot)$.
\end{lemma}

\section {Estimates   of the   norms  of the Birman-Schwinger  operator}
\setcounter{equation}{0}

 Let $H_0=-{\D}+x_1$ be the Stark operator.
We are  going to use  the  representation  of  $\exp(-itH_0)$  as  a product  of   different  factors,  one  of  which is  $\exp(it\D )$.
 One of   such formulas  was  discovered in  \cite{AH77} and is given by
\[
\lb{E2} e^{-itH_0}=e^{-itx_1}e^{it\D}e^{t^2{\pa \/\pa x_1} }e^{-i{t^3\/3}},\qqq
 \qq  \forall\  t\in \R.
\]
Another  representation  of $\exp(-itH_0)$  is
\[
 e^{-itH_0}=e^{-i{t^3\/12}} \bigl(e^{-itx_1/2}e^{it\D}e^{-itx_1/2}\bigr),\qqq
 \qq  \forall\  t\in \R.\notag
\]
One of the things   that make  this   formula  useful  is that the integral  kernel of the operator $e^{it
\D}$ on $L^2(\R^3)$ is  given explicitly
\[
\lb{eit} (e^{it\D})(x,y)={e^{-i3\pi/4}\/ (4\pi
t)^{3/2}}e^{i|x-y|^2/4t}, \ \ \ \ \ \ \ \ \ t> 0,
\]
$x,y\in \R^3$. The latter  observation allows  one to obtain a  nice representation for operators  that are related to
 $R_0(\l)=(H_0-\l)^{-1}$.  Indeed, let
\[\label{3.4}
{\frak  R}(\l,\z)=\int_0^\iy  e^{-it(H_0-\l)}t^{\z-1}dt.
\]  for all ${\rm Re}\,\z>0$.
If ${\rm Im }\,\l>0$, then the  integrals in   \eqref{3.4}  converge   (absolutely) in the  operator-norm topology. Moreover,
$$
R_0(\l)=i\,{\frak R}(\l,1).
$$
Throughout  the paper,  we  use the following  convenient   notation
$$
\Lambda=\l-2^{-1}(x_1+y_1).
$$
\begin{proposition} The operator $e^{-itH_0}$  is representable in the form
\[\lb{repr}
 e^{-itH_0}=e^{-i{t^3\/12}} \bigl(e^{-itx_1/2}e^{it\D}e^{-itx_1/2}\bigr),\qqq
 \qq  \forall\  t\in \R.
\]
 The integral kernel $r_\z(x,y,\l)$ of the operator $ {\frak  R}(\l,\z)$ equals
\[
\lb{E3}
r_\z(x,y,\l)={e^{-i{3\pi\/4}}\/\sqrt{(4\pi)^3}}\int_0^\iy
e^{{i\/4t}|x-y|^2}e^{-i{t^3\/12}}e^{it\,\Lambda}t^{\z-1}{dt\/t^{3/2}},\qquad {\rm Re}\,\z>3/2,
\]
for $x,y\in \R^3$, ${\rm Re}\,\z>3/2$ and $\l\in \C_+$.
\end{proposition}
{\it Proof}. The formula \eqref{repr} which implies  \eqref{E3}  can be proved by direct differentiation.
Indeed,  for any  $f\in C_0^\iy(\R^3)$,
\[\begin{aligned}
\frac{d}{dt}\Bigl(e^{-itx_1/2}e^{it\D}e^{-itx_1/2}\Bigr)f=&\\
-\frac{ix_1}2\Bigl(e^{-itx_1/2}e^{it\D}e^{-itx_1/2}\Bigr)f+i\Bigl(e^{-itx_1/2} \D e^{it\D}e^{-itx_1/2}\Bigr)f&-\frac{i}2\Bigl(e^{-itx_1/2}e^{it\D}{x_1} e^{-itx_1/2}\Bigr)f=\\
-iH_0\Bigl(e^{-itx_1/2}e^{it\D}e^{-itx_1/2}\Bigr)f+i\Bigl(\bigl[e^{-itx_1/2},\D \bigr]e^{it\D}e^{-itx_1/2}\Bigr)f&-\frac{i}2\Bigl(e^{-itx_1/2}\bigl[e^{it\D},{x_1}\bigr]e^{-itx_1/2}\Bigr)f.\notag
\end{aligned}
\]
It remains to note that
\[
\bigl[e^{-itx_1/2},\D \bigr]=e^{-itx_1/2}\bigl(it\frac{\partial}{\partial x_1}+t^2/4\bigr),\qquad \text{and}\qquad \bigl[e^{it\D},{x_1}\bigr]=2i t \frac{\partial}{\partial x_1}e^{it\D}.\notag
\]
$\BBox$

\bigskip

Let us now define the operators
$$
{\frak R}_1(\l,\z)=\int_0^1 e^{-it(H_0-\l)}t^{\z-1}dt\quad \text{and} \quad {\frak R}_2(\l,\z)=\int_1^\infty e^{-it(H_0-\l)}t^{\z-1}dt.
$$

\begin{proposition} \lb{interpolProp3.2}
The operators ${\frak R}_1(\l,\z)$  and ${\frak R}_2(\l,\z)$ are  bounded  if ${\rm \,Re}\, \z>0$ and ${\rm \,Re}\, \z<0$
correspondingly. Moreover,
\[
||{\frak R}_1(\l,\z)||\leq \frac1{{\rm \,Re}\, \z},\qquad {\rm \,Re}\, \z>0
\]
\[
||{\frak R}_2(\l,\z)||\leq \frac1{|{\rm \,Re}\, \z|},\qquad {\rm \,Re}\, \z<0
\]
The integral kernels  of the operators $
{\frak R}_1(\l,\z)\, \text{and} \, \,{\frak R}_2(\l,\z)
$ equal
\[\lb{r1r2}\begin{aligned}
\r_1(x,y;\l, \z)={e^{-i{3\pi\/4}}\/\sqrt{(4\pi)^3}}\int_0^1
e^{{i\/4t}|x-y|^2}e^{-i{t^3\/12}}e^{it\,\Lambda}t^{\z-1}{dt\/t^{3/2}},\quad \text{and} \\
\quad \r_2(x,y;\l, \z)={e^{-i{3\pi\/4}}\/\sqrt{(4\pi)^3}}\int_1^\iy
e^{{i\/4t}|x-y|^2}e^{-i{t^3\/12}}e^{it\,\Lambda}t^{\z-1}{dt\/t^{3/2}}\end{aligned}
\]
for $x,y\in \R^3$ and $\l\in \C_+$. There exists a  finite   $C_\z(p)>0$ depending only  on ${\rm Re}\, \z$and $p\geq1$  such that
\[\lb{L2r1r2}
\begin{aligned}
\int|\r_1(x,y;\l, \z)|^2d\l<C_\z(2),\quad \forall {\rm \,Re}\, \z>2\quad \text{and} \\
\quad \int|\r_2(x,y;\l, \z)|^pd\l<C_\z(p),\quad \forall  p<\frac2{2{\rm Re}\,\z-3}, \quad  3/2\leq  {\rm \,Re}\, \z<2.
\end{aligned}
\]
 There exists another  $\tilde C_\z>0$ depending only  on ${\rm Re}\, \z$ such that
\[\lb{Liyr1r2}\begin{aligned}
\sup_{x,y,\l}|\r_1(x,y;\l, \z)|<\tilde C_\z,\quad \forall {\rm \,Re}\, \z>3/2\quad \text{and}\\
\quad \sup_{x,y,\l}|\r_2(x,y;\l, \z)|<\tilde C_\z,\quad \forall {\rm \,Re}\, \z<3/2.\end{aligned}
\]
\end{proposition}
{\it Proof.}  All statements of this proposition are trivial. One only needs to  explain relations \eqref{L2r1r2}, which follow  from the  fact that, as  functions of $\l$, the kernels
$\r_1$  and $\r_2$  are  Fourier transforms  of   functions  that  could be estimated  by  $t^{\z-5/2}\chi(t)$, where $\chi $  is the characteristic  function of  either  $[0,1]$ or $[1,\iy)$.
In this sense, proving  the estimate involving   $\r_2$ is  more  difficult, because, additionally,  one  needs  to observe   that  $t^{\z-5/2}\chi(t)$  belongs to $L^q[1,\infty)$
with $q=p/(p-1)$. $\BBox$

\bigskip

The following result  helps  us to turn the information  provided   by \eqref{L2r1r2} into the information about the integral of  the norms  of
the Birman-Schwinger operators.

\begin{proposition}
Let $\eta(x,y,\l) $  be a measurable function on $\R^3\times \R^3\times \R$  such that
\[\lb{condeta}
||\eta||_{\iy,p}^p=\sup_{x,y}\int_\R |\eta(x,y,\l)|^pd\l<\iy,\qquad  p\geq 2.
\]
Let $T_\l$  be  the integral operator, whose kernel is $\eta(\cdot,\cdot,\l)$. Finally, let $W_1$ and $W_2$  be  two  functions from the space $L^2({\R^3})$.
Then  $W_1 T_\l W_2$ is a Hilbert-Schmidt operator  for almost every $\l\in \R$. Moreover,
\[
\int_\R ||W_1T_\l W_2||_{{\frak S}_2}^p \, d\l \leq ||\eta||_{\iy,p}^p||W_1||^p_{L^2}\cdot ||W_2||^p_{L^2}.
\]
\end{proposition}

{\it Proof.}  The case $p=2$ is obvious. Let us  assume now that $p>2.$
Then
\[\begin{aligned}
\int_\R\Bigl(\int_{\R^3}\int_{\R^3}|W_1(x)|^2|\eta(x,y,\l)|^2|W_2(y)|^2dxdy\Bigr)^{p/2}& d\l\leq \notag\\  \Bigl(\int_{\R^3}\int_{\R^3}|W_1(x)|^2|W_2(y)|^2dxdy\Bigr)^{p/2-1}
\int_\R\Bigl(\int_{\R^3}& \int_{\R^3}|W_1(x)|^2|\eta(x,y,\l)|^p|W_2(y)|^2dxdy\Bigr)d\l.
\end{aligned}
\]
The  statement of the proposition follows. $\BBox$

\bigskip
\begin{corollary}
Let $\ve\in (0,1/2)$. Then
\[\lb{infty}
||W_1{\frak R}_2(\l,\z)W_2||\leq \frac1{\ve}\,\,||W_1||_{L^\iy}\cdot ||W_2||_{L^\iy},\qquad {\rm \,Re}\, \z=-\ve,
\] and
\[\lb{intp}
\int_\R ||W_1{\frak R}_2(\l,\z)W_2||_{{\frak S}_2}^p \, d\l \leq C_\z(p)||W_1||^p_{L^2}\cdot ||W_2||^p_{L^2},\quad \forall  p<\frac2{1-2\ve}, \quad   {\rm \,Re}\, \z=2-\ve.
\]
\end{corollary}

Interpolating between the two cases  considered  in this corollary, we obtain the following  very important statement.
\begin{lemma}
Let $\ve\in(0\,1/2)$  and let
$\theta=(1+\ve)/2$.
\[\lb{interpol1}
\int_\R ||W_1{\frak R}_2(\l,1)W_2||_{{\frak S}_{2/\theta}}^{p/\theta} \, d\l \leq \ve^{p-p/\theta} C_{2-\ve}(p)||W_1||^{p/\theta}_{L^{2/\theta}}\cdot ||W_2||^{p/\theta}_{L^{2/\theta}},\quad \forall  p<\frac2{1-2\ve}.
\]
\end{lemma}

{\it Proof.}
Let us  take an arbitrary measurable operator-valued  function
$G(\cdot)$   such that
$$
||G||_{dual}:=\Bigl(\int_\R||G(\lambda)||_{ {\frak S}_{\frac{2}{2-\theta}}}^{\frac{p}{p-\theta}}d\l\Bigr)^{\frac{p-\theta}{p}}<\infty.
$$
For each $\l\in \R$,  the value $G(\l)$ is an operator in $L^2(\R^3)$.
Set  now $g(\l)=||G(\lambda)||_{ {\frak S}_{\frac{2}{2-\theta}}}$ and $Q(\l)=g(\l)^{-1}G(\l)$. Then the ${\frak S}_{\frac{2}{2-\theta}}$-norm of
$Q(\l)\in {\frak S}_{\frac{2}{2-\theta}}$ equals $1$  and    $g\in L^{\frac{p}{p-\theta}}(\R)$. Now we define
\[\begin{aligned}
f(\z)=
\int_{\R}\Tr\bigl[ |W_1|^{\frac{\z+\ve}{2\theta}}{\frak R}_2(\l,\z) |W_2|^{\frac{\z+\ve}{2\theta}} Q(\l)|Q(\l)|^{-1+\frac{4-\z-\ve}{4-2\theta}}\bigr] |g(\l)|^{\frac{2p-(\z+\ve)}{2(p-\theta)}}d\l.\notag
\end{aligned}
\]
Observe that,  for  any $t\in \R$,
\[
f(-\ve+it)=\int_{\R}\Tr\bigl[ |W_1|^{\frac{it}{2\theta}}{\frak R}_2(\l,-\ve+it) |W_2|^{\frac{it}{2\theta}} Q|Q|^{-1+\frac{4-it}{4-2\theta}}\bigr] |g(\l)|^{\frac{2p-it }{2(p-\theta)}}d\l.\notag
\]
Hence,   due to\eqref{infty},
\[\lb{line1}
|f(-\ve+it)|\leq \ve^{-1} ||g||^{\frac{p}{p-\theta}}_ { L_{\frac{p}{p-\theta}}}
\]
Similarly, since  the relation
\[
f(2-\ve+it)=\int_{\R}\Tr\bigl[ |W_1|^{\frac{2+it}{2\theta}}{\frak R}_2(\l,2-\ve+it) |W_2|^{\frac{2+it}{2\theta}} Q|Q|^{-1+\frac{2-it}{4-2\theta}}\bigr] |g(\l)|^{\frac{2p-2-it}{2(p-\theta)}}d\l.\notag
\]
implies the inequality
\[
|f(2-\ve+it)|\leq  ||g||^{\frac{p-1}{p-\theta}}_ { L_{\frac{p}{p-\theta}}}\Bigl( \int_{\R}||\,|W_1|^{\frac{2+it}{2\theta}}{\frak R}_2(\l,2-\ve+it) |W_2|^{\frac{2+it}{2\theta}}
||_{{\frak S}_2}^p d\l \Bigr)^{1/p},\notag
\]
we obtain from  \eqref{intp} that
\[\lb{line2}
|f(2-\ve+it)|\leq   C^{1/p}_{2-\ve}(p) ||g||^{\frac{p-1}{p-\theta}}_ { L_{\frac{p}{p-\theta}}}||W_1||^{1/\theta}_{L^{\frac{2}{\theta}}}||W_2||^{1/\theta}_{L^{\frac{2}{\theta}}}.
\]
It  follows  now from \eqref{line1} and \eqref{line2} by the three lines theorem, that
\[\lb{3lines}
|f(1)|\leq   \ve^{\theta-1} C^{\theta/p}_{2-\ve}(p) ||g||_ { L_{\frac{p}{p-\theta}}}||W_1||_{L^{\frac{2}{\theta}}}||W_2||_{L^{\frac{2}{\theta}}}.
\]
On the other  hand,
$$
f(1)=\int_{\R}\Tr\bigl[ |W_1|{\frak R}_2(\l,1) |W_2| G(\l)\bigr] d\l.
$$
Therefore, \eqref{3lines}
will turn into \eqref{interpol1}, once we take $G(\l)=w(\l) \Bigl| |W_1|{\frak R}_2(\l,1) |W_2| \Bigr|^{2/\theta-1}\Omega^*(\l)$
with  $$\Omega(\l)=|W_1|{\frak R}_2(\l,1) |W_2| \cdot \Bigl| |W_1|{\frak R}_2(\l,1) |W_2| \Bigr|^{-1}\,\,\text{and} $$
$$
\omega(\l)=|| |W_1|{\frak R}_2(\l,1) |W_2| ||_{{\frak S}_{2/\theta}}^{(p-2)/\theta}.
$$
$\BBox$

\bigskip

 Observe now, that  $p/\theta$ in \eqref{interpol1} is any number  satisfying
$$
\frac2{1+\ve}\leq \frac{p}\theta<\frac4{(1+\ve)(1-2\ve)}=\frac4{1-\ve-2\ve^2}.
$$
In particular, we  can choose $p/\theta=4/(1-\ve)$. Thus, we  obtain  the following
\begin{theorem} Let  $\ve\in (0,1/2)$, let $ p=\frac4{1-\ve}$ and let $q=\frac4{1+\ve}$. Then there  exists a   constant $C(\ve)>0$ such that
\[
\int_\R ||W_1{\frak R}_2(\l,1)W_2||_{{\frak S}_{q}}^{p} \, d\l \leq C(\ve)||W_1||^{p}_{L^q}\cdot ||W_2||^{p}_{L^q}.
\]
\end{theorem}

\bigskip

We  work  with the operator $W_1{\frak R}_1(\l,\z)W_2$ in the same  way. First, we formulate the  following  consequence of Proposition~\ref{interpolProp3.2}.

\begin{corollary}
Let $\ve\in (0,1/2)$. Then
\[
||W_1{\frak R}_1(\l,\z)W_2||\leq \frac1{\ve}\,\,||W_1||_{L^\iy}\cdot ||W_2||_{L^\iy},\qquad {\rm \,Re}\, \z=\ve,
\] and
\[
\int_\R ||W_1{\frak R}_1(\l,\z)W_2||_{{\frak S}_2}^2 \, d\l \leq C_{2+\ve}(2)||W_1||^2_{L^2}\cdot ||W_2||^2_{L^2}, \qquad   {\rm \,Re}\, \z=2+\ve.
\]
\end{corollary}

Interpolating between these two cases, we derive the estimate
\[
\int_\R ||W_1{\frak R}_1(\l,1)W_2||_{{\frak S}_{2/\theta}}^{2/\theta} \, d\l \leq C_\ve||W_1||^{2/\theta}_{L^{2/\theta}}\cdot ||W_2||^{2/\theta}_{L^{2/\theta}},
\]
where $\theta\in (0,1)$  is  the  number satisfying  the relation $\theta(2+\ve)+(1-\theta)\ve=1$. Put  differently
$\theta=(1-\ve)/2$.  Observe now, that
$$
 \frac{2}\theta=\frac4{1-\ve}.
$$
 Thus, we  obtain  the following  theorem.
\begin{theorem} Let  $\ve\in (0,1/2)$  and let $ p=\frac4{1-\ve}$. Then there  exists a   constant $C_\ve>0$ such that
\[
\int_\R ||W_1{\frak R}_1(\l,1)W_2||_{{\frak S}_{p}}^{p} \, d\l \leq C_\ve||W_1||^{p}_{L^p}\cdot ||W_2||^{p}_{L^p}.
\]
\end{theorem}

Finally, since $$R_0(\l)=i\bigl[{\frak R}_1(\l,1)+{\frak R}_2(\l,1)\bigr],$$
we conclude that the following assertion  can  be  made about    the Birman-Schwinger operators.
\begin{theorem}\lb{BSint} Let  $\tilde\ve\in (0,1/2)$, let $ p=\frac4{1-\tilde\ve}$ and let $q=\frac4{1+\tilde\ve}$. Then there  exists a   constant $C(\tilde\ve)>0$ such that
\[
\int_\R ||W_1R_0(\l+i\t)W_2||_{{\frak S}_{p}}^{p} \, d\l \leq C(\tilde\ve)\Bigl(||W_1||^{p}_{L^p}\cdot ||W_2||^{p}_{L^p}+||W_1||^{p}_{L^q}\cdot ||W_2||^{p}_{L^q}\Bigr),
\]
for any $\t\geq 0.$
\end{theorem}

\bigskip

Theorem~\ref{BSint} would not be so  useful in applications without  the following  result.

\begin{proposition}\label{zeros}
Let $a(\cdot)$ be an analytic function on $\C_+=\{\Im \l>0\}$  satisfying
\begin{equation}
\label{eq:ass1}
a(\l) = 1 + o(|\l|^{-1})
\qquad\text{as}\,\, |\l|\to\infty \, \text{in}\,\,\, \C_+.
\end{equation}
Assume that   there   is a family  of positive  functions $f_\ve\in L^1(\R)$,   $0<\ve<\varepsilon_0$,  such that
\[\label{a<f}\ln|a(\l+i\varepsilon )|\leq  f_\ve(\l),\qquad \forall \l\in\R,\qquad  \forall \ve\in(0,\ve_0).\]
Then the zeros $\l_j$ of $a(\cdot)$ in $\C_+$, repeated according to multiplicities, satisfy
\begin{equation}
\label{eq:zeros}
\sum_j \Im \l_j \leq \frac{1}{2\pi}\sup_{0<\varepsilon<\varepsilon_0} \int_\R f_\ve(\l)\,d\l.
\end{equation}
\end{proposition}

{\it Proof.} Consider   first  the  function $a_\varepsilon(\l)=a(\l+i\varepsilon)$ for $0<\ve<\ve_0$.  Note that  zeros  of  $a_\varepsilon$  are  the points $\l_j-i\varepsilon$. Consequently,  the Blaschke product  for $a_\varepsilon$ in $\C_+$ is
$$
B_\ve(\l) = \prod_{\Im\l_j>\ve} \frac{\l+i\ve-\l_j}{\l-i\ve-\overline{\l_j}} \,.
$$
Since  $a_\ve(\l)/B_\ve(\l)$ is  analytic and non-zero in $\C_+=\{\Im \l>0\}$, the    function  $\log (a_\ve(\l)/B_\ve(\l))$ exists and is analytic there. For $R>0$ we denote by $C_R$ the contour which consists of the interval $[-R,R]$, traversed from left to right, and the circular part $\Gamma_R:=\{\l\in {\mathbb C}:\, |\l|=R \,, \Im \l >0\}$, traversed counterclockwise.
Then
$$
\int_{C_R} \log \frac{a_\ve(\l)}{B_\ve(\l)}\,d\l = 0 \,,
$$
and, therefore,
\begin{equation}
\label{eq:zerosproof1}
\Re \int_{-R}^{R} \log \frac{a_\ve(x)}{B_\ve(x)}\,dx + \Re \int_{\Gamma_R} \log\frac{a_\ve(\l)}{B_\ve(\l)}\,d\l = 0 \,.
\end{equation}
We note that $|B_\ve(x)|=1$ if $x\in\R$ and, therefore,
\[\label{eq:zerosproof2}
\begin{aligned}
\Re \int_{-R}^{R} \log \frac{a_\ve(x )}{B_\ve(x )}\,dx
& = \int_{-R}^{R} \ln \left| \frac{a_\ve(x )}{B_\ve(x )} \right| \,dx \\
& = \int_{-R}^{R} \ln \left| a_\ve(x ) \right| \,dx \,.
\end{aligned}\]
(We denote by $\ln$ the natural logarithm to distinguish it from the particular branch of the complex logarithm $\log$ chosen before.) On the other hand, by \eqref{eq:ass1} and $B_\ve(\l) = 1+O(|\l|^{-1})$ (note that the zeros $\l_j$ are contained in a bounded set as a consequence of \eqref{eq:ass1}), both $\log a_\ve(\l)$ and $\log B_\ve(\l)$ are well-defined for all sufficiently large $|\l|$ and we have, for all sufficiently large $R$,
\begin{equation}
\label{eq:zerosproof3}
\Re \int_{\Gamma_R} \log\frac{a_\ve(\l)}{B_\ve(\l)}\, d\l = \Re \int_{\Gamma_R} \log a_\ve(\l)\, d\l - \Re \int_{\Gamma_R} \log B_\ve(\l)\,d\l \,.
\end{equation}
We conclude from \eqref{eq:zerosproof1}, \eqref{eq:zerosproof2} and \eqref{eq:zerosproof3} that
\begin{equation}
\label{eq:traceformula}
\Re \int_{\Gamma_R} \log B_\ve (\l)\, d\l = \int_{-R}^{R} \ln \left| a_\ve(x ) \right| \,dx + \Re \int_{\Gamma_R} \log a_\ve (\l)\, d\l
\end{equation}
for all sufficiently large $R$.
We assume that   $| \l_j-i\ve|<R$ for all $j$.  Since
$$
\log B_\ve (\l) = 2i \sum_{\Im\l_j>\ve} \frac{\ve- \Im  \l_j}{\lambda} + O( (\lambda)^{-2}) \,,
$$
we get
\begin{equation}
\begin{aligned}
\label{eq:tf1}
\int_{\Gamma_R} \log B_\ve (\l) \,d\l & = \\
& = -2\pi \sum_{\Im\l_j>\ve} \left(\ve - \Im  \l_j\right) + O( R^{-1})
\qquad\text{as}\ R\to\infty .
\end{aligned}
\end{equation}
On the other hand, by \eqref{eq:ass1},
\[
\begin{aligned}
\label{eq:tf2}
\Re \int_{\Gamma_R} \log a_\ve (\l)\,d\l = o(1)
\qquad\text{as}\ R\to\infty \,.
\end{aligned}
\]
Moreover, by \eqref{a<f},
\begin{equation}
\begin{aligned}
\label{eq:tf3}
\int_{-R}^{R} \ln \left| a_\ve(x ) \right| \,dx & \leq \int_{-R}^{R} f_\ve(\l)d\l \leq  \int_{-\infty}^{\infty} f_\ve(\l)d\l\,.
\end{aligned}
\end{equation}
Relations \eqref{eq:traceformula},  \eqref{eq:tf1}, \eqref{eq:tf2} and \eqref{eq:tf3} imply
\begin{equation}
\label{eq:zeros2}
\sum_j ( \Im \l_j -\ve)_+\leq \frac{1}{2\pi} \int_\R f_\ve(\l)\,d\l\leq \frac{1}{2\pi}\sup_{0<\ve<\ve_0} \int_\R f_\ve(\l)\,d\l.
\end{equation}
Inequality \eqref{eq:zeros} now follows from \eqref{eq:zeros2} by  the monotone convergence theorem. $\,\,\,\,\,\,\,\,\,\,\,\,\,\,\,\BBox$

\section{Proof of Theorem~\ref{main}}
It is  sufficient to  consider  the case $V\in C_0^\infty(\R^3)$.
We
apply   Proposition~\ref{zeros} with $$a(\lambda)=\det_5(I+W_1R_0(\l)W_2)$$  and  $f_\ve(\l)=C||W_1R_0(\l+i\ve)W_2||_{{\mathfrak S}_p}^p$ where $C$ is the constant  from \eqref{det<p}  with $n=5$ and $4<p\leq5$. Note  that  the  zeros of $a(\l)$
are  eigenvalues  of $H.$

We  use Theorem~\ref{BSint} to conclude that
%%%%%%%%%%%%%%%%%%%%%%%%%%%%%%%%%%%%%%%%%%%%%%%%%%%%%%%%%%%%%%%%%%%%%%%%%%%%%%%%%
$$
\int_\R f_\ve(\l)d\l\leq C
\Bigl[\Bigl(\int_{\R^3}|V|^{p/2}dx\Bigr)^{2}+\Bigl(\int_{\R^3}|V|^{q/2}dx\Bigr)^{2p/q}\Bigr].
$$
where  $q=4/(1+\tilde\ve)$ and  $\tilde\ve\in(0,1/5]$ is  such that $p=4/(1-\tilde\ve)$.
Also, it is  established in the last section   that $Y_0(\l)$ satisfies the  following estimate.
\begin{theorem}
Let ${\rm Im}\l\geq0$. Then  for any $p>11$,  there exists  a  positive constant  $C_p>0$ depending only on $p$  such that
\[
||Y_0(\l)||_{{\mathfrak S}_4}\leq\frac{ C_p}{1+|\l|^{1/4}}\Bigl(\int_{{\mathbb R}^3}(1+|x|)^p|V|^2 dx\Bigr)^{1/2}.
\]
\end{theorem}
It  is  clear from this theorem  that $a(\l)=1+O(|\l|^{-5/4})$, as $|\l|\to \infty$.
So, all conditions  of Proposition~\ref{zeros} are  fulfilled, and therefore, Theorem~\ref{main} follows.  $\BBox$

\section{Non-integral bounds  for the  Birman-Schwinger operator}

Another  consequence of Proposition~\ref{interpolProp3.2} is the following  statement
\begin{proposition} Let $\ve>0.$
Let $W_1$  and $W_2$ be two functions on $\R^3$.
Then
\[
\begin{aligned}
||W_1{\frak R}_1(\l,\z)W_2||\leq \frac1{\ve}||W_1||_{L^\iy}||W_2||_{L^\iy},\qquad {\rm Re}\,\z=\ve\\
||W_1{\frak R}_2(\l,\z)W_2||\leq \frac1{\ve}||W_1||_{L^\iy}||W_2||_{L^\iy},\qquad {\rm Re}\,\z=-\ve.
\end{aligned}
\]
Moreover, there  is a constant $C_\ve>0$ such that
\[
\begin{aligned}
||W_1{\frak R}_1(\l,\z)W_2||_{{\frak S}_2}\leq  C_{\ve}||W_1||_{L^2}   ||W_2||_{L^2},\qquad {\rm Re}\,\z=3/2+\ve, \\
||W_1{\frak R}_2(\l,\z)W_2||_{{\frak S}_2}\leq  C_{\ve}||W_1||_{L^2}    ||W_2||_{L^2},\qquad {\rm Re}\,\z=3/2-\ve.
\end{aligned}
\]
\end{proposition}

\bigskip

The standard interpolation ( that  has  been used  already in this paper)  leads  to
\begin{theorem} Let $\ve\in(0,1/2).$
Let $W_1$  and $W_2$ be two functions on $\R^3$.
Then
\[
\begin{aligned}
||W_1{\frak R}_1(\l,1)W_2||_{{\frak S}_{2/\theta}}\leq \tilde C_{\ve}||W_1||_{L^{2/\theta}}   ||W_2||_{L^{2/\theta}},\qquad \theta=2(1-\ve)/3,\\
||W_1{\frak R}_2(\l,1)W_2||_{{\frak S}_{2/\tilde\theta}}\leq\tilde C_{\ve}||W_1||_{L^{2/\tilde\theta}}   ||W_2||_{L^{2/\tilde\theta}},\qquad \tilde \theta=2(1+\ve)/3.
\end{aligned}
\]
\end{theorem}
\bigskip

Finally, applying the triangle inequality, we obtain
\begin{corollary}
Let $\ve\in(0,1/2).$
Let $W_1$  and $W_2$ be two functions on $\R^3$.
Then there  exists a  positive constant $C(\ve)>0,$ such that
\[\lb{interpol5.40}
\begin{aligned}
||W_1{ R}_0(\l)W_2||_{{\frak S}_p}\leq  C({\ve})\bigl( ||W_1||_{L^p}  ||W_2||_{L^p}+||W_1||_{L^q}   ||W_2||_{L^q}\bigr),\\ \qquad \text{where}
\qquad p=3/(1-\ve),
\quad q=3/(1+\ve).
\end{aligned}
\]
\end{corollary}

\section{Interpolation for ${\rm Im}\,\l>0.$}

Our starting point is the formula
\[\label{s3.4}
{\frak  R}(\l,\z)=\int_0^\iy  e^{-it(H_0-\l)}t^{\z-1}dt.
\]  for all ${\rm Re}\,\z>0$.
If ${\rm Im }\,\l>0$, then the  integrals in   \eqref{s3.4}  converge   (absolutely) in the  operator-norm topology. Moreover,
$$
R_0(\l)=i\,{\frak R}(\l,1).
$$
We remind the reader that
the integral kernel $r_\z(x,y,\l)$ of the operator $ {\frak  R}(\l,\z)$ equals
\[
\lb{sE3}
r_\z(x,y,\l)={e^{-i{3\pi\/4}}\/\sqrt{(4\pi)^3}}\int_0^\iy
e^{{i\/4t}|x-y|^2}e^{-i{t^3\/12}}e^{it\,\Lambda}t^{\z-1}{dt\/t^{3/2}},
\]
where $$
\Lambda=\l-2^{-1}(x_1+y_1),\qquad  x,y\in \R^3,\quad \text{and}\quad \l\in \C_+.$$

\begin{proposition} \lb{sinterpolProp3.2}
The operator ${\frak R}(\l,\z)$  is  bounded  if ${\rm \,Re}\, \z>0$
and ${\rm \,Im}\, \l>0$. Moreover,  there  exists a positive constant $C_{\z}>0$ depending only   on ${\rm \,Re}\, \z$, such that
\[
||{\frak R}(\l,\z)||\leq \frac{C_{ \z}}{|{\rm \,Im}\, \l|^{{\rm \,Re}\, \z}},\qquad {\rm \,Re}\, \z>0.
\]
There exists a  finite   $C_\z>0$ depending only  on ${\rm Re}\, \z$  such that
\[\lb{sL2r1r2}
\begin{aligned}
\Bigl(\int_\R|r_\z(x,y,\l+i\g)|^2d\l\Bigr)^{1/2}<\frac{C_{ \z}}{\g^{{\rm \,Re}\, \z-2}},\qquad \forall {\rm \,Re}\, \z>2.
\end{aligned}
\]
 There exists another  $\tilde C_\z>0$ depending only  on ${\rm Re}\, \z$ such that
\[\lb{sLiyr1r2}
\sup_{x,y,\l}|r_\z(x,y,\l)|<\frac{\tilde  C_{ \z}}{|{\rm \,Im}\, \l|^{{\rm \,Re}\, \z-3/2}},\qquad \forall {\rm \,Re}\, \z>3/2.
\]
\end{proposition}

\bigskip

Let us  now turn this information into the information about the Birman-Schwinger operators.

\begin{corollary}
Let $\ve>0$  and $\t\geq 2+\ve$. Then
\[\lb{sinfty}
||W_1{\frak R}(\l,\z)W_2||\leq\frac{C_{ \ve}}{|{\rm \,Im}\, \l|^{\ve}}\,\,||W_1||_{L^\iy}\cdot ||W_2||_{L^\iy},\qquad {\rm \,Re}\, \z=\ve,
\] and
\[\lb{sintp}
\Bigl(\int_\R ||W_1{\frak R}(\l+i\g,\z)W_2||_{{\frak S}_2}^2 \, d\l \Bigr)^{1/2}\leq \frac{C_{  \t}}{\g^{\t-2}}||W_1||_{L^2}\cdot ||W_2||_{L^2},\quad \forall   {\rm \,Re}\, \z= \t.
\]
\end{corollary}

Interpolating between the two cases  considered  in this corollary, we obtain the following  very important statement.
\begin{lemma}\label{slemma}
Let $0<\ve<1$,  $\t\geq 2+\ve$  and let
$\theta=(1-\ve)/( \t-\ve)$. Then
\[\lb{sinterpol1}
\int_\R ||W_1{\frak R}(\l+i\g,1)W_2||_{{\frak S}_{2/\theta}}^{2/\theta} \, d\l \leq  \frac{C_\ve^{2/\theta-2} C^{2}_{\t}}{\g^{2/\theta-4}}||W_1||^{2/\theta}_{L^{2/\theta}}\cdot ||W_2||^{2/\theta}_{L^{2/\theta}}.
\]
\end{lemma}

{\it Proof.}
Let us  take an arbitrary measurable operator-valued  function
$Q(\cdot)$   such that
$$
||Q||_{dual}:=\Bigl(\int_\R||Q(\lambda)||_{ {\frak S}_{\frac{2}{2-\theta}}}^{\frac{2}{2-\theta}}d\l\Bigr)^{\frac{2-\theta}{2}}<\infty.
$$
For each $\l\in \R$,  the value $Q(\l)$ is an operator in $L^2(\R^3)$.
Now we define
\[\begin{aligned}
f(\z)=
\int_{\R}\Tr\bigl[ |W_1|^{\frac{\z-\ve}{1-\ve}} {\frak R}(\l+i\g,\z) |W_2|^{\frac{\z-\ve}{1-\ve}} Q(\l)|Q(\l)|^{-1+\frac{2 \t-\z-\ve}{( \t-\ve)(2-\theta)}}\bigr]d\l.\notag
\end{aligned}
\]
Observe that,  for  any $t\in \R$,
\[
f(\ve+it)=\int_{\R}\Tr\bigl[ |W_1|^{\frac{it}{1-\ve}} {\frak R}(\l+i\g,\ve+it) |W_2|^{\frac{it}{1-\ve}} Q(\l)|Q(\l)|^{-1+\frac{2 (\t-\ve)-it}{( \t-\ve)(2-\theta)}}\bigr]d\l.\notag
\]
Hence,   due to  \eqref{sinfty},
\[\lb{sline1}
|f(\ve+it)|\leq \frac{C_{ \ve}}{\g^{\ve}}\, ||Q||_{dual}^{\frac{2}{2-\theta}}
\]
Similarly, since  the relation
\[
f(\t+it)=\int_{\R}\Tr\bigl[ |W_1|^{\frac{\t-\ve+it}{1-\ve}} {\frak R}(\l+i\g,\t+it) |W_2|^{\frac{\t-\ve+it}{1-\ve}} Q(\l)|Q(\l)|^{-1+\frac{\t-\ve-it}{( \t-\ve)(2-\theta)}}\bigr] d\l.\notag
\]
implies the inequality
\[
|f(\t+it)|\leq  ||Q||^{\frac{1}{2-\theta}}_ {dual}\Bigl( \int_{\R}||\,|W_1|^{\frac{\t-\ve+it}{1-\ve}} {\frak R}(\l+i\g,\t+it) |W_2|^{\frac{\t-\ve+it}{1-\ve}}
||_{{\frak S}_2}^2 d\l \Bigr)^{1/2},\notag
\]
we obtain from  \eqref{sintp} that
\[\lb{sline2}
|f(\t+it)|\leq   \frac{C_{  \t}}{\g^{\t-2}} ||Q||^{\frac{1}{2-\theta}}_ {dual}||W_1||^{1/\theta}_{L^{\frac{2}{\theta}}}||W_2||^{1/\theta}_{L^{\frac{2}{\theta}}}.
\]
It  follows  now from \eqref{line1} and \eqref{line2} by the three lines theorem, that
\[\lb{s3lines}
|f(1)|\leq   \frac{C_\ve^{1-\theta} C^{\theta}_{\t}}{\g^{1-2\theta}}||Q||_ {dual}||W_1||_{L^{\frac{2}{\theta}}}||W_2||_{L^{\frac{2}{\theta}}}.
\]
On the other  hand,
$$
f(1)=\int_{\R}\Tr\bigl[ |W_1| {\frak R}(\l+i\g,1) |W_2| Q(\l)\bigr] d\l.
$$
Therefore, \eqref{s3lines}
will turn into \eqref{sinterpol1}, once we take $Q(\l)= \Bigl| |W_1| {\frak R}(\l+i\g,1) |W_2| \Bigr|^{2/\theta-1}\Omega^*(\l)$
with  $$\Omega(\l)=|W_1| {\frak R}(\l+i\g,1) |W_2| \cdot \Bigl| |W_1| {\frak R}(\l+i\g,1) |W_2| \Bigr|^{-1}. $$
$\BBox$

\bigskip

It is more  convenient  to formulate  Lemma~\ref{slemma}
in the following  way.
\begin{theorem}\lb{sthm6.4}
Let $p>4$ and $\g>0$. Then there  exists a constant $C_p>0$ depending only on $p$ such that
\[\lb{sinterpol1}
\int_\R ||W_1 R_0(\l+i\g)W_2||_{{\frak S}_p}^{p} \, d\l \leq  \frac{C_p}{\g^{p-4}}||W_1||^{p}_{L^p}\cdot ||W_2||^{p}_{L^p}.
\]
\end{theorem}

Another  consequence of Proposition~\ref{sinterpolProp3.2} is the following  statement
\begin{proposition} Let $\ve>0$ and let $\t\geq  3/2+\ve$.
Let $W_1$  and $W_2$ be two functions on $\R^3$.
Then there  exists a  constant $C_\ve>0$  such that
\[
\begin{aligned}
||W_1{\frak R}(\l,\z)W_2||\leq \frac{C_\ve}{|{\rm Im}\l|^\ve}||W_1||_{L^\iy}||W_2||_{L^\iy},\qquad {\rm Re}\,\z=\ve
\end{aligned}
\]
Moreover, there  is a constant $C_\t>0$ such that
\[
\begin{aligned}
||W_1{\frak R}_1(\l,\z)W_2||_{{\frak S}_2}\leq  \frac{C_\t}{|{\rm Im}\l|^{\t-3/2}}||W_1||_{L^2}   ||W_2||_{L^2},\qquad {\rm Re}\,\z=\t.
\end{aligned}
\]
\end{proposition}

The standard interpolation ( that  has  been used  already in this paper)  leads  to
\begin{lemma} Let $\ve\in(0,1)$, $\t\geq 3/2+\ve$  and $\theta=(1-\ve)/(\t-\ve)$.
Let $W_1$  and $W_2$ be two functions on $\R^3$.
Then
\[
\begin{aligned}
||W_1{\frak R}(\l,1)W_2||_{{\frak S}_{2/\theta}}\leq  \frac{C_\ve^{1-\theta}C^\theta_\t}{|{\rm Im}\l|^{1-3\theta/2}}||W_1||_{L^{2/\theta}}   ||W_2||_{L^{2/\theta}}.
\end{aligned}
\]
\end{lemma}

Put  differently, we  obtain
\begin{theorem} \lb{p3}
Let $p>3.$
Let $W_1$  and $W_2$ be two functions on $\R^3$.
Then there  exists a  positive constant $C_p>0,$ such that
\[\lb{interpol5.40}
\begin{aligned}
||W_1{ R}_0(\l)W_2||^p_{{\frak S}_p}\leq  \frac{C_p}{|{\rm Im }\l|^{p-3}} ||W_1||^p_{L^p}  ||W_2||^p_{L^p}
\end{aligned}
\]
\end{theorem}

\section{Proof of Theorem~\ref{p>1}}

In this   section,  we   establish   some bounds  on  the sums of  the powers of  imaginary parts  of the eigenvalues  of the operator $H$.
  First,  we  prove   the  following  statement which could be  viewed as a  generalization  of  Theorem~\ref{main}.
\begin{theorem}\lb{lift} Let $p>4$.
Let $V\in L^{p/2}(\R^3)$ be   a  bounded   complex-valued potential.
Then  there  exists  a   positive constant $C_p>0$ depending only on $p$,  such that   for any $\g>0$, the  eigenvalues $\lambda_j$   of the operator $H$   satisfy  the  estimate
\[
\lb{c8}
 \sum_j
 ({\rm Im}\, \lambda_j-\g)_+\leq \frac{C_p}{\g^{p-4}}\Bigl[\int_{{\mathbb R}^3}|V(x)|^{p/2}\,dx\Bigr]^2,
\]
\end{theorem}

{\it Proof.} It is  sufficient to consider  the case $V\in C_0^\infty(\R^3)$.
We
apply   Proposition~\ref{zeros} with
$$a(\lambda)=\det_n\Bigl(I+Y_0(\l+i\g)\Bigr),\qquad  n-1\leq p\leq n,$$  and  $f_\ve(\l)=C||Y_0(\l+i(\gamma+\ve))||_{{\mathfrak S}_p}^p$ where $C$ is the constant  from \eqref{det<p}. Note  that  a point $\l\in {\mathbb C}_+$ is a  zero of $a(\l)$
if and  only  if  $\l+i\g$ is an  eigenvalue  of $H.$
Note also that  Theorem~\ref{sthm6.4}  implies the  inequality
%%%%%%%%%%%%%%%%%%%%%%%%%%%%%%%%%%%%%%%%%%%%%%%%%%%%%%%%%%%%%%%%%%%%%%%%%%%%%%%%%
$$
\int_\R f_\ve(\l)d\l\leq  \frac{C_p}{\g^{p-4}}\Bigl[\int_{{\mathbb R}^3}|V(x)|^{p/2}\,dx\Bigr]^2.
$$
It  is also clear  that $a(\l)=1+O(|\l|^{-5/4})$, as $|\l|\to \infty$. So, all conditions  of Proposition~\ref{zeros} are fulfilled.  Therefore, Theorem~\ref{lift} follows.
$\BBox$

\bigskip

{\it Proof of Theorem~\ref{p>1}}.   It is enough to consider  the case $V\in C^\infty_0(\R^3)$. First, we observe that  according to Theorem~\ref{p3}, there  exists a  positive constant  $c_p>0$  depending only on $p$ such that
 $|{\rm Im}\l_j|<c_p\Bigl[\int_{{\mathbb R}^3}|V(x)|^{p/2}\,dx\Bigr]^{2/(p-3)} $  for all $j$.
Now we multiply \eqref{c8}  by $\g^{q-2}$  and  integrate  the  resulting  inequality   with respect to $\g$  from $0$  to $\g_0=c_p\Bigl[\int_{{\mathbb R}^3}|V(x)|^{p/2}\,dx\Bigr]^{2/(p-3)} $
\begin{equation}\lb{dlinnaya}\begin{aligned}
 \sum_j
 \int_0^\infty ({\rm Im}\, \lambda_j-\g )_+\g^{q-2} d\g\leq \Bigl[\int_{{\mathbb R}^3}|V(x)|^{p/2}\,dx\Bigr]^2 \int_0^{\g_0} \frac{C_p}{\g^{p-4}}\g^{q-2}d\g,
\end{aligned}\end{equation}
A simple change of the variable in the corresponding  integral leads to  the equality
\[\lb{st64}
 \int_0^\infty ({\rm Im}\, \lambda_j-\g )_+\g^{q-2} d\g=|{\rm Im}\, \lambda_j|^q\int_0^\infty (1-\g )_+\g^{q-2} d\g.
\]
$\BBox$

\section{Resolvent operator. Revised}

 Let $H_0=-{\D}+x_1$ be the free Stark operator.
The representation
\[
 e^{-itH_0}=e^{-i{t^3\/12}}\Bigl(e^{-itx_1/2}e^{it\D}e^{-itx_1/2 }\Bigr),\qqq
 \qq  \forall\  t\in \R,\notag
\]
also   implies that  the closure of $e^{zx_1}e^{-zH_0}$  is  representable  as described  below.

\begin{proposition} Let ${\rm Re}\,\,z\geq0$. If $K(z)$  is  the closure of $e^{zx_1}e^{-zH_0}$, then
\[
\lb{AH78}\bigl( K(z) f, g\bigr)=e^{{z^3\/12}}\Bigl(\bigl(e^{z\D}e^{-zx_1/2 }\bigr)f,e^{\bar zx_1/2} g\Bigr),\qqq \forall f,g\in C_0^\iy(\R^3),
 \qq  \forall\   {\rm Re}\,\,z\geq0.
\]
\end{proposition}

Note   that  \eqref{AH78}  could be formally written in the  following  (dubious) form
\[\lb{dubious}
 e^{-zH_0}=e^{{z^3\/12}}\Bigl(e^{-zx_1/2}e^{z\D}e^{-zx_1/2 }\Bigr),\qqq
 \qq  \forall\ {\rm Re}\,\,z\geq0.
\]
{\it Proof  }.
The formula \eqref{AH78}  follows  from the observation that  the  quantity
$\bigl( e^{-zH_0}E_{H_0}(a,b)f,g\bigr)
$  depends  analytically on $z$  for any $-\infty<a<b<\infty$, $f\in L^2({\mathbb R}^3)$ and $g\in C^\infty_0({\mathbb R}^3)$.
 On the other hand,  for the same $a$,  $b$, $f$  and $g$, the quantity
$$
\bigl(e^{z\D}e^{-z^2{\pa \/\pa x_1} }e^{{z^3\/3}}E_{H_0}(a,b)f, e^{-\bar zx_1}g\bigr)
$$
depends  analytically   on $z$ in  the right  half-plane $\{z:\,\,\text{Re}\,z>0\}.$  Due to the  fact  that it is  also continuous  up to the boundary  of the half-plane, we obtain  from \eqref{E2}  that
$$
\bigl( e^{-zH_0}E_{H_0}(a,b)f,g\bigr)=\bigl(e^{z\D}e^{-z^2{\pa \/\pa x_1} }e^{{z^3\/3}}E_{H_0}(a,b)f,e^{-\bar zx_1} g\bigr),\qquad  {\rm Re}\,z\geq 0.
$$
One can drop  the spectral projection $E_{H_0}(a,b)$, provided  that the resulting relation will be  written in the form
$$
\bigl( K(z) f, e^{-\bar zx_1}g\bigr)=\bigl(e^{z\D}e^{-z^2{\pa \/\pa x_1} }e^{{z^3\/3}}f,e^{-\bar zx_1} g\bigr),\qquad  {\rm Re}\,z\geq 0.
$$
The latter   relation implies \eqref{AH78}, because
$$e^{{z^3\/3}}\bigl(e^{z\D}e^{-z^2{\pa \/\pa x_1} }f,e^{-\bar zx_1} g\bigr)=e^{{z^3\/12}}\bigl(e^{z\D}e^{-z x_1/2}f,e^{-\bar zx_1/2} g\bigr)\qquad \forall\, f,g\in C^\infty_0({\mathbb R}^3).$$
$\BBox$

\bigskip

As  we see, $ e^{-tH_0}$  is not a  continuous   operator. However,  (the closure of) the  product $e^{z x_1}e^{-z H_0}$   could be  viewed  as  a bounded operator   for all ${\rm Re}\,z\geq0$, due to the fact that $\Delta$ is a negative operator.

 Observe now  that  for  any $-\infty<a<b<\infty$,  the  product of  resolvent  operator
$R_0(\l)=(H_0-\l)^{-1}$  and  the spectral projection $E_{H_0}(a,b)$ can be  written as  the sum of  two  integrals
\[\lb{6.29}
 R_0(\l)E_{H_0}(a,b)=\int_0^{1}  e^{-t(H_0-\l)}E_{H_0}(a,b)dt+ i\int_0^{\infty}  e^{-i(t-i)(H_0-\l)}E_{H_0}(a,b)dt .
\]
While the first integral   converges   for all $\l$,
the  second  integral in  the right  hand  side of   \eqref{6.29}  converges (absolutely) in the operator-norm topology only  for ${\rm Im \l}>0.$
We  will often drop the  projection $E_{H_0}(a,b)$  and  write  formally that
\[\lb{6.29b}
 R_0(\l)=\int_0^{1}  e^{-t(H_0-\l)}dt+ i\int_0^{\infty}  e^{-i(t-i)(H_0-\l)}dt .
\]
Now  we  are going to obtain a  useful  representation for the  integral kernels of  the   operators  in the right hand   side of \eqref{6.29b}.
In particular, we will show that  the  quadratic   forms of these operators are  well-defined on $C_0^\infty(\R^3)$.

We will also  show that  the  two  terms  in the right hand  side  of \eqref{6.29b}   could  be  viewed  as the   values of   two   families   of operators
$
{\frak T}_0(\l,\z)$ and ${\frak B}(\l,\z)
$
at $\z=1.$ These  families  will depend on $\z$  analytically, which will allow   us  to  interpolate.

Let us deal with the first  term in the right hand side  of  \eqref{6.29b}.
%%%%%%%%%%%%%%%%%%%%%%%%%%%%%%
In order to  use  an interpolation,  we introduce  the family of (unbounded) operators  $${\frak T}_0(\l,\z)=e^{\z^2}\int_0^{1}  e^{-t(H_0-\l)}t^{\z-1}dt,$$   depending    analytically on the parameter $\z$  after being multiplied   by a  spectral projection $E_{H_0}(a,b)$.
Observe that, at least  formally,
$${\frak T}_0(\l,\z)=e^{\z^2}\int_0^{1} e^{t^3/12} e^{-tx_1/2}e^{t\D}e^{-tx_1/2} e^{t\l}t^{\z-1}dt.$$
Set now
 $${\frak T}(\l,\z)=(1+|x_1|)^{-(\z-\ve)/(\t-\ve)}{\frak T}_0(\l,\z)(1+|x_1|)^{-(\z-\ve)/(\t-\ve)},\qquad \ve>0,\,\,\t>5/2,$$
and define ${\mathcal P}(x)=1+e^{-x_1/2}$.
It is  easy to see that if $W_1$ and $W_2$ are two functions on $\R^3$, and ${\rm Re}\,\z=\ve>0$, then
\[\lb{i1}
||W_1{\frak T}(\l,\z)W_2||\leq \frac{(1+e^{{\rm Re}\,\l})}\ve||{\mathcal P} W_1||_{L^\iy} \cdot ||{\mathcal P} W_2||_{L^\iy} .
\]
Let  us  denote the  kernel of the operator ${\frak T}_0(\l,\z)$  by $\n(x,y;\l,\z)$.
Recall  again that the integral  kernel of the operator   $e^{t
\D}$  on $L^2(\R^3)$ is the   function
\[
 (e^{t\D})(x,y)={1\/ (4\pi
t)^{3/2}}e^{-|x-y|^2/4t}, \ \ \ \ \ \ \ \ \ t> 0,
\]
where $x,y\in \R^3$.
Obviously,
the function
\[\lb{definitionnu}
\n(x,y,\l,z)={e^{\z^2}\/\sqrt{(4\pi)^3}}\int_0^1
e^{{-1\/4t}|x-y|^2}e^{{t^3\/12}-2^{-1}t(x_1+y_1))} e^{\l t}\,{t^{\z-1}dt\/t^{3/2}},
\end{equation}
satisfies the inequality
\[\lb{obviousnu<}
|\n(x,y;\l,\z)|\leq C_{{\rm Re}\,\z}(1+e^{{\rm Re}\l}){\mathcal P}(x){\mathcal P}(y)
\quad
\text{ for\,\, all}  \quad {\rm Re}\,\z>3/2.
\]
 It  turns out  that $\n$   decays  as $|\l|\to\iy$  in the half-plane $\{\l:\,\,\,{\rm Re}\,  \l<\alpha\}$  for any $\a\in {\mathbb R}$. In order to obtain an estimate that  shows   such a  behavior of $\nu$,
we   write $\n$
 in the  form
\begin{equation}
\lb{nu=}\begin{aligned}
\n(x,y,\l,z)={e^{\z^2}\/\l\sqrt{(4\pi)^3}}\int_0^1
e^{{-1\/4t}|x-y|^2}e^{{t^3\/12}-2^{-1}t(x_1+y_1))}\bigl(\frac{d e^{\l t}}{dt}\bigr){t^{\z-1}dt\/t^{3/2}},
\end{aligned}
\end{equation}
for $x,y\in \R^3$ and $\l\in \C_+\setminus \{0\}$.
Integrating by parts in \eqref{nu=}, we obtain
\begin{equation}
\lb{nparts}\begin{aligned}
\n(x,y;\l,\z)=\,\,\,\,\,\,\,\,\,\,\,\,\,\,\,\,\,\,\,\,\,\,\,\,\,\,\,\,\,\,\,\,\,\,\,\,\,\,\,\,\,\,\,\,\,\,\,\,\,\,\,\,\,\,\,\,\,\,\,\,\,\,\,\,\,\,\,\,\,\,\,\,\,\,\,\,\,\,\,\,\,\,\,\,\,\,\,\,\,\,\,\,\,\,\,\,\,\,\,\,\,\,\,\,\,\,\,\,\,\,\,\,\,\,\,\,\,\,\,\,\,\,\,\,\,\,\,\,\,\,\,\,\,\,\,\,\,\,\,\,\,\,\,\,\,\,\,\,\,\,\,\,\,\,\,\,\,\,\,\,\,\,\,\,\,\,\,\,\,\,\,\,\\
{e^{\z^2}\/\l\sqrt{(4\pi)^3}}\int_0^1
e^{{-1\/4t}|x-y|^2+{t^3\/12}+t\,\Lambda}\Bigl(2^{-1}(x_1+y_1)+\Bigl(\frac5{2}-\z\Bigr)t^{-1}-\frac{|x-y|^2+t^4}{4t^2}   \Bigr){t^{\z-1}  dt\/t^{3/2}}+\\+{1\/\l\sqrt{(4\pi)^3}}
e^{{-1\/4}|x-y|^2+{1\/12}+\,\Lambda}
\end{aligned}
\end{equation}
for $x,y\in \R^3$ and $\l\in \C_+$.    The formula \eqref{nparts}  (combined with the  inequality \eqref{obviousnu<})  leads  to  the  estimate
\begin{lemma}\label{ocenkanu} Let ${\mathcal P}_1(x)=(1+|x_1|)(1+e^{-x_1/2})$.
There  exists a  positive  constant $C_\t>0$  such that
\[\lb{ocenkas2}
|\n(x,y;\l,\z)|\leq C_\t\frac{(1+e^{{\rm Re}\l})}{1+| \l|}{\mathcal P}_1(x){\mathcal P}_1(y),\qquad  {\rm Re}\,\z=\t>5/2
\]
\end{lemma}

\begin{corollary} The Hilbert-Schmidt  norm of the operator $W_1{\frak T}(\l,\z)W_2$  satisfies the estimate
\[\lb{i2}
||W_1{\frak T}(\l,\z)W_2||_{{\frak S}_2}\leq C_\t\frac{(1+e^{{\rm Re}\,\l})}{1+|\l|}||{\mathcal P} W_1||_{L^2} \cdot ||{\mathcal P} W_2||_{L^2} ,\qquad  {\rm Re}\,\z=\t>5/2.
\]
\end{corollary}
Interpolating between  \eqref{i1}  and \eqref{i2} we obtain

\begin{proposition}
Let $p>5$ and let ${\mathcal P}(x)=1+e^{-x_1/2}$.  Then the  ${\frak S}_p-$  norm of the operator $W_1{\frak T}(\l,1)W_2$  satisfies the estimate
\[\lb{Tp}
||W_1{\frak T}(\l,1)W_2||_{{\frak S}_p}^p\leq C_p\frac{(1+e^{{\rm Re}\,\l})^p}{(1+|\l|)^2}||{\mathcal P} W_1||_{L^p}^p \cdot ||{\mathcal P} W_2||_{L^p}^p .
\]
\end{proposition}

{\it Proof.} We  use   \eqref{i1} with $0<\ve<1$ and \eqref{i2} with $\t>5/2$. The  previously used interpolation technique leads to
\eqref{Tp}  with   $p=2/\theta$ where  $\theta\in(0,1)$  satisfies the relation $\ve(1-\theta)+\t\theta=1$. Put  differently,
$\theta=(1-\ve)/(\t-\ve)$, which implies  that $p$  can be any number  greater than $5$. $\BBox$
%%%%%%%%%%%%%%%%%%%%%%%%%%%%%%%%%%%%%%%%%%%%%%%

\medskip

This proposition immediately implies

\begin{proposition}   Let $V$ be a complex-valued  function on $\R^3$ and
let $p>5$.  Assume that $W_1$ and $W_2$ satisfy the relations $W_1=|V|^{1/2} $ and
$V=W_1 W_2$.  Then the ${\frak S}_p-$  norm of the operator $W_1{\frak T}_0(\l,1)W_2$  satisfies the estimate
\[\lb{fT0Lp}
||W_1{\frak T}_0(\l,1)W_2||_{{\frak S}_p}^p\leq C_p\frac{(1+e^{{\rm Re}\,\l})^p}{(1+|\l|)^2}\Bigl(\int_{\R^3} (1+e^{-x_1/2})^p(1+|x_1|)^2 |V|^{p/2}dx\Bigr)^2 .
\]
\end{proposition}

Similarly, one can deal with the second term in the right hand side of \eqref{6.29b}.
First of all we note that
$$ i\int_0^{\infty}  e^{-i(t-i)(H_0-\l)}dt= i\int_0^{\infty} e^{-i{(t-i)^3\/12}} e^{-(1+it)x_1/2} e^{(1+it)\D}e^{-(1+it)x_1/2}e^{(1+it)\l}dt.$$
Following  the steps of  our  work  with the first  term in \eqref{6.29b}, we  introduce the operators
$$
{\frak B}(\l,\z)= i\int_0^{\infty}  e^{-i(t-i)(H_0-\l)}t^{\z-1}dt
$$
for  ${\rm Re}\,\z>0.$  First, we  observe that the following statement  holds  true.
\begin{proposition} Let $\ve\in (0,1)$.
Let $W_1$  and $W_2$  be two functions on $\R^3$. Let $\psi (x)=e^{-x_1/2}$.  Then
\[\lb{B<e}
||W_1{\frak B}(\l,\z)W_2||\leq C_\ve\frac{ e^{{\rm Re}\,\l+(({\rm Im}\,\l)_-+1)^2}}{\bigl(1+({\rm Im}\,\l)_+\bigr)^\ve} ||\psi  W_1||_{L^\iy}||\psi  W_2||_{L^\iy},\qquad  \text{for \, all}\quad {\rm Re}\,\z=\ve.
\]
\end{proposition}

{\it Proof.} One only needs to estimate the integral
$$
 \int_0^{\infty}\bigl| e^{-i{(t-i)^3\/12}}e^{(1+it)\l} t^{\z-1}\bigr|dt= e^{{\rm Re}\,\l}\int_0^{\infty}e^{{-3t^2+1\/12}}e^{-t{\rm Im}\,\l} t^{\ve-1}dt.
$$
Assume that ${\rm Im}\,\l>0$. Then
$$
 e^{{\rm Re}\,\l}\int_0^{\infty}e^{{-3t^2+1\/12}}e^{-t{\rm Im}\,\l} t^{\ve-1}dt\leq C e^{{\rm Re}\,\l}\int_0^{\infty}e^{-t{\rm Im}\,\l} t^{\ve-1}dt\leq
\frac{ C_\ve e^{{\rm Re}\,\l}}{({\rm Im}\,\l)_+^\ve}.
$$
If ${\rm Im}\,\l<0$, then
$$
 e^{{\rm Re}\,\l}\int_0^{\infty}e^{{-3t^2+1\/12}}e^{-t{\rm Im}\,\l} t^{\ve-1}dt\leq  e^{{\rm Re}\,\l}\int_0^{1}e^{{-3t^2+1\/12}}e^{-t{\rm Im}\,\l} t^{\ve-1}dt
+ e^{{\rm Re}\,\l}\int_1^{\infty}e^{{-3t^2+1\/12}}e^{-t{\rm Im}\,\l} dt.
$$
Both integrals in the right  hand side   can be  estimated  in a very  simple way:
 $$\int_0^{1}e^{{-3t^2+1\/12}}e^{-t{\rm Im}\,\l} t^{\ve-1}dt\leq C_\ve e^{({\rm Im}\,\l)_-},$$
and
$$
\int_1^{\infty}e^{{-3t^2+1\/12}}e^{-t{\rm Im}\,\l} dt\leq \int_{-\infty}^{\infty}e^{{-3t^2+1\/12}}e^{-t{\rm Im}\,\l} dt\leq C e^{({\rm Im}\,\l)^2}.
$$
$\BBox$

\medskip

 One  can also  provide a proof  of   the  following  statement.
\begin{proposition}
 The  integral kernel   of the operator $ {\frak B}(\l,\z)= i\int_0^{\infty}  e^{-i(t-i)(H_0-\l)}t^{\z-1}dt$  is the function
$$
\eta(x,y;\l,\z)={e^{-i{\pi\/4}}\/\sqrt{(4\pi)^3}}\int_0^\iy
e^{{it-1\/4(t^2+1)}|x-y|^2}e^{-i{(t-i)^3\/12}}e^{i(t-i)\,\Lambda}{t^{\z-1}dt\/(t-i)^{3/2}},
$$
where the agreement  about  the  choice of the  branch  of $(t-i)^{3/2}$  is that $(t-i)^{3/2}\bigr|_{t=0}=e^{-i3\pi/4}$.
\end{proposition}

{\it Proof.}
  We use \eqref{dubious}  to  obtain
\[
\eta(x,y;\l,\z)=(2\pi)^{-3}i\int_{0}^\infty e^{i(t-i)\Lambda -\frac{i(t-i)^3}{12}}\Bigl(\int_{{\mathbb R}^3} e^{i\tilde p(x-y)}e^{-i(t-i)|\tilde p|^2}d\tilde p\Bigr)\,  t^{\z-1} dt.\notag
\]
Now, the  statement of the proposition  follows from the  fact that
\[
(2\pi)^{-3}\int_{{\mathbb R}^3} e^{i\tilde p(x-y)}e^{-i(t-i)|\tilde p|^2}d\tilde p= {e^{-i3\pi/4}\/\sqrt{(4\pi)^3}}\frac{
e^{{i\/4(t-i)}|x-y|^2}}{(t-i)^{3/2}}.\notag
\]
$\BBox$

\bigskip

Note that
$$
|\eta(x,y;\l,\z)|\leq {1\/\sqrt{(4\pi)^3}}\int_0^\iy
e^{{(1-3t^2)\/12}}\bigl|e^{(1+it)\,\Lambda}\bigr|\,\,t^{{\rm Re}\,\z-1}dt={e^{{\rm Re}\,\Lambda}\/\sqrt{(4\pi)^3}}\int_0^\iy
e^{{(1-3t^2)\/12}}e^{-t{\rm Im }\,\Lambda}\,t^{{\rm Re}\,\z-1}dt.
$$
Consequently,  we  state the following result.

\begin{proposition} Let  $\t={\rm Re}\,\z>3/2$  and let $\psi(x)=e^{-x_1/2}$.  Then
\[\lb{eta<}
|\eta(x,y;\l,\z)|\leq C_\t{e^{{\rm Re}\,\lambda+2({\rm Im}\l)_-^2}\/(1+({\rm Im}\l)_+)^\t}\p(x)\p(y).
\]
\end{proposition}

{\it Proof.} Let  ${\rm Im}\,\l<0$. Then  by the  Schwarz inequality,
$$
|\eta(x,y;\l,\z)|\leq C_\t e^{{\rm Re}\,\Lambda}\Bigl(\int_0^\iy
e^{{(1-3t^2)\/12}}e^{-2t{\rm Im }\,\lambda}\, dt\Bigr)^{1/2}\leq C_\t e^{{\rm Re}\,\Lambda+2({\rm Im}\l)_-^2}.
$$
If ${\rm Im}\,\l>0$, then
$$
|\eta(x,y;\l,\z)|\leq {e^{{\rm Re}\,\Lambda }e^{{1\/12}}\/\sqrt{(4\pi)^3}}\int_0^\iy
e^{-t{\rm Im }\,\lambda}\,t^{\t-1}dt\leq C_\t{  e^{{\rm Re}\, \Lambda }\/  ({\rm Im }\,\l  )^\t }.
$$
$\BBox$

\bigskip

\begin{corollary}
 Let  $\t={\rm Re}\,\z>3/2$  and let $\psi(x)=e^{-x_1/2}$.  Then
\[\lb{fBL2}
||W_1{\frak B}(\l,\z)W_2||_{{\frak S}_2}\leq C_\t{e^{{\rm Re}\,\lambda+2({\rm Im}\l)_-^2}\/(1+({\rm Im}\l)_+)^\t}||\p W_1||_{L^2}||\p W_2||_{L^2}
\]
\end{corollary}

Interpolating   between \eqref{B<e}  and  \eqref{fBL2}, we derive
\begin{proposition} Let $p>3$
and let $\psi(x)=e^{-x_1/2}$.  Then  there  exists  a constant  depending only  on $p$  such that
\[\lb{fBLp}
||W_1{\frak B}(\l,1)W_2||_{{\frak S}_p}\leq C_p{e^{{\rm Re}\,\lambda+2(({\rm Im}\l)_-+1)^2}\/(1+({\rm Im}\l)_+)}||\p W_1||_{L^p}||\p W_2||_{L^p}
\]
\end{proposition}

Finally, combining \eqref{fT0Lp}  and \eqref{fBLp}, we obtain  by the triangle inequality  that  the   ${\frak S}_p$-norm of  $Y_0(\l)$
could be estimated as follows:
\begin{theorem}\lb{thm5.6}    Let $V$ be a complex-valued  function on $\R^3$ and
let $p>5$.  Assume that $W_1$ and $W_2$ satisfy the relations $W_1=|V|^{1/2} $ and
$V=W_1 W_2$.  Then the ${\frak S}_p-$  norm of the operator $Y_0(\l)$  satisfies the estimate
\[\lb{Y0Lp}
||Y_0(\l)||_{{\frak S}_p}^p\leq C_p\Bigl[\frac{1+e^{p{\rm Re}\,\l}}{(1+|\l|)^2}+{e^{p{\rm Re}\,\lambda+2p({\rm Im}\l)_-^2}\/(1+({\rm Im}\l)_+)^p}\Bigr]\Bigl(\int_{\R^3} (1+e^{-x_1/2})^p(1+|x_1|)^2 |V|^{p/2}dx\Bigr)^2 .
\]
\end{theorem}

\section{Jensen's  inequality for a  function analytic  in a  corner and  its  applications}

Here we prove the following  result about zeros of  an analytic  function.
\begin{proposition}\lb{PrN} Let $\varepsilon>0$ and $\a>0$  be two positive numbers. Let $a(z)$ be an analytic  function on  the  domain
$\Omega= \{z\in {\mathbb C} : \,\,{\rm Re}\, z \leq \a+\varepsilon, \,\, {\rm Im}\, z\geq  -\varepsilon\}$, having the asymptotics
$a(z) = 1 + o( |z|^{-2})$ as $|z|\to\infty$  in $\Omega$.  Assume also that
\[\lb{cond1ineq}
\ln |a(z)|\leq  \Bigl(\frac{(1+e^{p{\rm Re}z})}{(1+ |z|)^{2}}+\frac{ e^{2p({\rm Im} z)_-^2+p{\rm Re}z}}{(1+({\rm Im}z)_+)^p}\Bigr)^{1+\d}\cdot M,\qquad  \text{ if}\quad z\in \Omega,
\]
for  some $M>0,$ ,  $p>5$ and $\delta>0$,  which are  independent of $z$. Then the number $ N $ of zeros  of $a(z)$  in the domain $\{z\in {\mathbb C} : \,\,{\rm Re}\, z \leq \a, \,\, {\rm Im}\, z\geq  0\}$  satisfies
\[\lb{ocenkaN}
N\leq  \varepsilon^{-2} C_{p,\d}\cdot  M \Bigl[  e^{(1+\d)p(\a+\varepsilon)}\Bigl( \frac{\alpha+\varepsilon}{\varepsilon^{1+2\d}}+ (1+\ve^2)e^{2(1+\d)p\varepsilon^2} \Bigr)\Bigr],
\]
where $C_{p,\d}>0$ is independent of $\a$, $\varepsilon$  and $M.$
\end{proposition}

{\it  Proof}. The function   $\log(a(z)) $   it is not analytic in $\Omega$, due to the  possibility of  having zeros of $a(z)$  in $\Omega$. To  get rid of the zeros,  we introduce the following Blaschke product:
$$B(z) =\prod_ j \frac{(z -\a+ (i-1)\varepsilon)^2 -(z_j -\a+ (i-1 )\varepsilon)^2}{ (z -\a+ (i-1 )\varepsilon)^2 -(\bar z_j -\a- (1+i)\varepsilon)^2} ,$$
where $z_j$ are  zeros  of $a(z).$ It is   easy to see that  the  function
$
\log [a(z)/B(z)]
$  is analytic on $\Omega$, because  $B(z)$  vanishes  exactly at the points  $z=z_j$.  On the other  hand, $|B(z)|=1$  for all $z$ that  belong to the boundary of $\Omega$.

 Let $ C_R =\{z\in \Omega  :\,\, |z-\a| = R\}$, let $ I_R $ be  the interval $\{z\in \Omega:\,-\sqrt{R^2-\varepsilon^2}\leq {\rm Re}\,z-\a\leq \varepsilon  ,\, \,\text{and} \,\, {\rm Im}z = -\varepsilon \}$,  and let $J_R$ be  the interval $\{z\in \Omega: {\rm Re}\, z =\a+\varepsilon,\, \, \text{and}\,  -\varepsilon \leq  {\rm Im}\,z\leq  \sqrt{R^2-\varepsilon^2} \}$.
Define  $\Gamma_R=C_R\cup  I_R \cup J_R$ as a  traversed  counterclockwise  contour.
 Then $$\int_{\Gamma_R}\log [a(z) /B(z)] (z-\a +(i-1)\varepsilon)dz=0. $$
Consequently,
\begin{equation*}
\begin{aligned}\lim_{R\to\infty}{\rm Re}\int_{C_R}\log [B(z)/a(z)] (z-\a +(i-1)\varepsilon)dz=\lim_{R\to\infty}\int_{I_R}\log |a(z) | (z-\a +(i-1)\varepsilon)dz+\\
\lim_{R\to\infty}\int_{J_R}\log  |a(z) |(z-\a +(i-1)\varepsilon)dz, \end{aligned}
\end{equation*}
which implies that
\begin{equation}\lb{7.52}
\begin{aligned}\lim_{R\to\infty}{\rm Re}\int_{C_R}\log [B(z)] (z-\a +(i-1)\varepsilon)dz=\lim_{R\to\infty}\int_{I_R}\log |a(z) | (z-\a +(i-1)\varepsilon)dz+\\
\lim_{R\to\infty}\int_{J_R}\log  |a(z) |(z-\a +(i-1)\varepsilon)dz. \end{aligned}
\end{equation}
On the other hand,  due to the   expansion
$$
\log [B(z)] (z-\a +(i-1)\varepsilon)=-\frac{2i}z\sum_j{\rm Im }\, (z_j-\a+(i-1)\varepsilon)^2 +O(1/|z|^2),\qquad \text{as}\quad |z|\to\infty,
$$
the limit of the  integral in the left hand  side can be  easily computed. Namely,
$$
\lim_{R\to\infty}{\rm Re}\int_{C_R}\log [B(z)] (z-\a +(i-1)\varepsilon)dz=\pi\sum_j{\rm Im }\, (z_j-\a+(i-1)\varepsilon)^2.
$$
Therefore,
$$
\lim_{R\to\infty}{\rm Re}\int_{C_R}\log [B(z)] (z-\a +(i-1)\varepsilon)dz=2\pi\sum_j({\rm Im }\,z_j+\varepsilon)({\rm Re}z_j-(\a+\varepsilon))\leq -2\pi \varepsilon^2 N.
$$
Taking into account   the condition \eqref{cond1ineq},
we  obtain  from \eqref{7.52} that
\[\lb{2pie}\begin{aligned}
2\pi \varepsilon^2 N\leq M \int_{-\infty}^{\a+\varepsilon}  \Bigl(\frac{1+e^{pt}}{ (1+|t-i\varepsilon|)^{2}}+ e^{2p\varepsilon^2}e^{pt}\Bigr)^{1+\d} |t-\a-\varepsilon |dt+&\\ M \int_{-\varepsilon}^\infty
\Bigl( \frac{ 1+e^{p(\a+\ve)}}{(1+ |\alpha+\ve +it|)^{2}}&+\frac{e^{p(\a+\ve)}e^{2pt_-^2}}{(1+t_+)^p}\Bigr)^{1+\d}(t+\varepsilon)dt.
\end{aligned}
\]
Note that
\[\begin{aligned}
\int_{-\varepsilon}^\infty
\Bigl( \frac{ 1+e^{p(\a+\ve)}}{(1+ |\alpha+\ve +it|)^{2}}+\frac{e^{p(\a+\ve)}e^{2pt_-^2}}{(1+t_+)^p}\Bigr)^{1+\d}&(t+\varepsilon)dt
 \leq \\
 C_{\d,p}\Bigl[\frac{ 1+e^{(1+\d)p(\a+\ve)}}{(\a+\varepsilon)^{2\d}}&+ e^{(1+\d)p(\a+\ve)}(1+\ve+e^{2p(1+\d)\ve^2}\ve^2)\Bigr]
\end{aligned}\notag
\]
and
\[\begin{aligned}
\int_{-\infty}^{\a+\varepsilon}  \Bigl(\frac{1+e^{pt}}{ (1+|t-i\varepsilon|)^{2}}+ e^{2p\varepsilon^2}e^{pt}\Bigr)^{1+\d} |t-\a-\varepsilon |dt
\leq \\
C_{p,\d} \Bigl[\frac{(1+e^{(1+\d)p(\a+\varepsilon)})(\alpha+\varepsilon)}{\varepsilon^{1+2\d}}+ e^{2(1+\d)p\varepsilon^2}e^{(1+\d)p(\a+\varepsilon)}\Bigr].\notag
\end{aligned}
\]
Consequently, \eqref{2pie}  can be written in the form
$$
2\pi \varepsilon^2 N\leq  C_{\d,p}\cdot  M \Bigl[  e^{(1+\d)p(\a+\varepsilon)}\Bigl( \frac{\alpha+\varepsilon}{\varepsilon^{1+2\d}}+ (1+\ve^2)e^{2(1+\d)p\varepsilon^2} \Bigr)\Bigr].
$$
The proof is  completed. $\,\,\,\,\,\,\,\,\,\,\,\,\,\,\,\,\,\,\,\,\,\,\,\,\,\BBox$

\bigskip

We now can  apply  this proposition   to the function
$$
a(z)=\det_{n}(I+Y_0(z)).
$$
where $n$  is  such that $n-1\leq p(1+\d)\leq n$
Let's  remind the reader that
 according to   Theorem~\ref{thm5.6}   combined with the inequality \eqref{det<p}, there exists  a   positive  constant $C_{p,\d}>0$ depending on $p$ and $\d$  such that
\eqref{cond1ineq}   holds  with
$$
M=C_{p,\d}
\Bigl(\int_{\R^3} (1+e^{-x_1/2})^p(1+|x_1|)^2 |V|^{p/2}dx\Bigr)^{2(1+\d)}.
$$
Thus, Theorem~\ref{main2} follows  from Proposition~\ref{PrN}.

\section{ Individual   eigenvalue  bounds. Proof  of Theorem~\ref{individual}}

Here  we   obtain an estimate of the Hilbert-Schmidt  norm of the Birman-Schwinger operator that allows  us  to say something about the location of  eigenvalues of $H$ in the complex plane.

Let $H_0=-{\D}+x_1$ be the  free  Stark operator.
We are  going to use  the  representation  of  $\exp(-itH_0)$  as  a product  of   different  factors,  one  of  which is  $\exp(it\D )$.
Namely,
\[
\lb{individualE2} e^{-itH_0}=e^{-i{t^3\/12}} \Bigl(e^{-itx_1/2}e^{it\D}e^{-itx_1/2}  \Bigr),\qqq
 \qq  \forall\  t\in \R.\notag
\]
 On the other hand, the  resolvent  operator
$R_0(\l)=(H_0-\l)^{-1}$  can be     written as  the integral
\[
 R_0(\l)=i\int_0^\iy  e^{-it(H_0-\l)}dt.\notag
\]
If ${\rm Im }\,\l>0$, then this  integral  converges   (absolutely) in the  operator-norm topology.
We remind the reader that
$$
\Lambda=\l-2^{-1}(x_1+y_1).
$$
\begin{proposition} The integral kernel $r_0(x,y,\l)$ of the operator $R_0(\l)$ equals
\[
\lb{10E3}
r_0(x,y,\l)={e^{i\sqrt{ \,\Lambda}\,|x-y|}\/4\pi|x-y|}+\\{e^{-i{\pi\/4}}\/\sqrt{(4\pi)^3}}\int_0^\iy
e^{{i\/4t}|x-y|^2}\Bigl(e^{-i{t^3\/12}}-1\Bigr)e^{it\,\Lambda}{dt\/t^{3/2}},
\]
for $x,y\in \R^3$ and $\l\in \C_+$.
\end{proposition}
{\it Proof}.  Indeed, since
\[
 R_0(\l)=i\int_0^\iy e^{-i{t^3\/12}} \Bigl(e^{-itx_1/2}e^{it\D}e^{-itx_1/2}  \Bigr) e^{it\l}dt,\notag
\]
we come to the conclusion that
\[\lb{R087}
 R_0(\l)=i\int_0^\iy  e^{-itx_1/2}e^{it\D}e^{-itx_1/2} e^{it\l}dt + i\int_0^\iy (e^{-i{t^3\/12}}-1) \Bigl(e^{-itx_1/2}e^{it\D}e^{-itx_1/2}  \Bigr) e^{it\l}dt.
\]
It remains to observe that  the two terms  in \eqref{10E3} are the integral kernels of the operators in the right hand side of \eqref{R087}. $\BBox$

\bigskip

  The following  estimate  plays  a key role in the  arguments  of this section.
\begin{lemma}\label{l2.1} Let $W_1$ and $W_2$  be  two  functions from the space $L^3({\R^3})$.
Then
\[\label{CW}
\frac1{(4\pi)^2}\int_{\R^3}\int_{\R^3}\frac{|W_1(x)|^2|W_2(y)|^2}{|x-y|^2}dxdy\leq C^2 ||W_1||_{L^3}^{2}||W_2||_{L^3}^{2}.
\]
where $C>0$ is independent of $W_1$ and $W_2$.
\end{lemma}

{ \it Proof.} Note  that  the function
$
{W_1(x)W_2(y)}/(4\pi |x-y|)
$
is the integral  kernel of the operator $T=W_1(-\Delta)^{-1}W_2.$ According to the Cwikel-Lieb-Rozenblum  inequality (see \cite{Cw},\cite{Lieb} and \cite{Roz}),
the number $n(s, T)$ of singular  values  of $T$
 lying   to the right of $s>0$ satisfies the  relation
$$
n(s,T)\leq C s^{-3/2}||W_1||_{L^3}^{3/2}||W_2||_{L^3}^{3/2}
$$
with a constant $C>0$ independent of $W_1$  and $W_2$. In particular, it implies the bound  $||T||\leq  C^{2/3} ||W_1||_{L^3}||W_2||_{L^3}$.
It remains to note that
$$
\frac1{(4\pi)^2}\int_{\R^3}\int_{\R^3}\frac{|W_1(x)|^2|W_2(y)|^2}{|x-y|^2}dxdy=||T||^2_{{\mathfrak S}_2}=2\int_0^{||T||}  n(s, T) s  ds.
$$
{\BBox}

\medskip

\begin{corollary} \lb{corollary3.3} Let $W_1$ and $W_2$  be  two  functions from the space $L^3({\R^3})$.
Then
$$
||W_1(-\Delta-\l)^{-1}W_2||_{{\mathfrak S}_2}^2\leq  C^2 ||W_1||_{L^3}^{2}||W_2||_{L^3}^{2},\qquad  \forall \l\in \overline{\mathbb C}_+,
$$
where $C$ is the same as in \eqref{CW}.
\end{corollary}

{\it Proof.}
The function $e^{i\sqrt{\l}|x-y|}/(4\pi|x-y|)$  is the kernel of  the operator $(-\Delta-\l)^{-1}$. Consequently,
$$
||W_1(-\Delta-\l)^{-1}W_2||_{{\mathfrak S}_2}^2\leq \frac1{(4\pi)^2}\int_{\R^3}\int_{\R^3}\frac{|W_1(x)|^2|W_2(y)|^2}{|x-y|^2}dxdy.
$$
$\BBox$

\medskip

Let us  now introduce   the  following  convenient  notations
$$
\mu(\lambda,x,y)={e^{-i{\pi\/4}}\/\sqrt{(4\pi)^3}}\int_0^\iy
e^{{i\/4t}|x-y|^2}\Bigl(e^{-i{t^3\/12}}-1\Bigr)e^{it\,\Lambda}{dt\/t^{3/2}},
$$
and
$$
\mu_1(\lambda,x,y)=\mu_0(\Lambda,x,y)-\mu_0(\lambda,x,y),
$$
where
$$
\mu_0(\lambda,x,y):={1\/4\pi|x-y|}e^{i\sqrt \l|x-y|}.
$$
In these  notations,
$$
r_0(x,y,\l)=\mu_0(\lambda,x,y)+\mu_1(\lambda,x,y)+\mu(\lambda,x,y).
$$
This representation of the  integral kernel leads  to the  corresponding  decomposition of  the  resolvent   operator
$$
R_0(\l)={\frak F}_0(\lambda)+{\frak F}_1(\lambda)+{\frak F}(\lambda).
$$
We  also  need to introduce the characteristic  function   $ \chi_{\lambda}(x) $ of  the  set $\{x\in \R^3:\,\,\,|x_1|<|\l|/2\}$.
\begin{lemma}\label{l2.2} Let $V$, $W_1$ and $W_3$  be three   functions  on ${\mathbb R}^3$ such that
$
V=W_2W_1,$  and $ |W_1|=|W_2|.
$
Let ${\frak F}_1(\lambda)$  be the integral operator on $L^2({\mathbb R}^3)$ with the kernel $\mu_1(\lambda,x,y)$.  Let $C_0$ be the best  constant in \eqref{CW}. Then
$$
||W_1{\frak F}_1(\lambda)W_2-\chi_{\lambda}W_1{\frak F}_1(\lambda)W_2\chi_{\lambda}||_{\mathfrak S_2}\leq \frac{2^{9/4}C_0}{(2+|\lambda|)^{1/4}}\Bigl(\int_{\R^3}(1+|x_1|)^{3/4}|V|^{3/2}dx\Bigr)^{2/3}.
$$
\end{lemma}
{\it Proof.} Note  that  due to the   fact  that  square of   the Hilbert-Schmidt norm of  an operator is the   integral  of the  square  of  its  kernel, we    have
$$
||W_1{\frak F}_1(\lambda)W_2||_{\mathfrak S_2}\leq  2||W_1{\frak F}_0(0)W_2||_{\mathfrak S_2}\leq  2C_0||W_1||_{L^3}||W_2||_{L^3}.
$$
 The statement of the lemma  follows  from   the  simple  fact  that
$$
||W_1(1-\chi_\lambda)||_{L^3}+||W_2(1-\chi_\lambda)||_{L^3}\leq \frac{2^{5/4}}{(2+|\lambda|)^{1/4}}\Bigl(\int_{\R^3}(1+|x_1|)^{3/4}|V|^{3/2}dx\Bigr)^{1/3}.
$$
The proof is completed. $\,\,\,\,\,\,\,\,\,\,\,\,\,\,\,\,\BBox$

\bigskip

\begin{lemma}\label{l2.3}
Let ${\frak F}_1(\lambda)$  be the operator  with the kernel $\mu_1(\lambda,x,y)$.  Let $\l \in \overline\C_+\setminus \{0\} $. Then
$$
||\chi_{\lambda}W_1{\frak F}_1(\lambda)W_2\chi_{\lambda}||_{\mathfrak S_2}\leq \frac{1}{4\pi(1+|\l|)^{1/4}}\Bigl(\int_{\R^3}(1+|x_1|)^{3/2}|V|dx\Bigr).
$$
\end{lemma}

{\it Proof.}  It is  easy to see  that the  kernel $\mu_1$  satisfies the   estimate
\[\lb{3.7}
|\chi_\lambda(x)\mu_1(\l,x,y)\chi_\lambda(y)|\leq \frac{1}{4\pi(1+|\l|)^{1/4}}\bigl((1+|x_1|)(1+|y_1|)\bigr)^{3/4}.
\]
Indeed,  since
$$
\mu_1(\l,x,y)={e^{i\sqrt{ \,\Lambda}\,|x-y|}\/4\pi|x-y|}-{e^{i\sqrt{ \l}\,|x-y|}\/4\pi|x-y|},\qquad \Lambda=\l-(x_1+y_1)/2,
$$
we obtain that
$$
|\chi_\lambda(x)\mu_1(\l,x,y))\chi_\lambda(y)|\leq {\Bigl|\sqrt{ \,\Lambda}-\sqrt{ \l}\Bigr|\cdot|x-y|\/4\pi|x-y|}\leq \frac{\sqrt2}{16\pi|\lambda|^{1/2}}|x_1+y_1|,
$$
for $(|x_1|+|y_1|)<|\l|$. The latter implies that
$$
|\chi_\lambda(x)\mu_1(\l,x,y))\chi_\lambda(y)|\leq \frac1{8\pi},\qquad \text{for}\quad |\l|\leq1.
$$ That proves \eqref{3.7}  for $|\l|\leq1$.  The  estimate \eqref{3.7} in the case $|\l|>1$  follows  from   the inequality
$$
 \frac{\sqrt2}{16\pi|\lambda|^{1/2}}|x_1+y_1|\leq  \frac{1}{4\pi(1+|\lambda|)^{1/4}}|x_1+y_1|^{3/4},\qquad  \text{for}\,\,|\l|>1.
$$
  The statement of  the lemma  immediately  follows  from \eqref{3.7}  and the definition of the Hilbert-Schmidt  norm. $\,\,\,\,\,\,\,\,\,\,\,\,\,\BBox$

\medskip
\begin{corollary}\label{col3.5}
Let ${\frak F}_1(\lambda)$  be the operator  with the kernel $\mu_1(\lambda,x,y)$.  Let $C_0$ be the best  constant in \eqref{CW}.  Let $\l \in \overline\C_+\setminus \{0\} $. Then
$$
||W_1{\frak F}_1(\lambda)W_2||_{\mathfrak S_2}\leq \frac{2^{9/4}C_0}{(2+|\lambda|)^{1/4}}\Bigl(\int_{\R^3}(1+|x_1|)^{3/4}|V|^{3/2}dx\Bigr)^{2/3}+ \frac{1}{4\pi(1+|\l|)^{1/4}}\Bigl(\int_{\R^3}(1+|x_1|)^{3/2}|V|dx\Bigr).
$$
\end{corollary}

\medskip

Let us now  consider the operator $W_1{\frak F}_0(\lambda)W_2$, where  ${\frak F}_0(\lambda)$ is the operator with   the integral kernel $\mu_0(\l,x,y)$.  According to  Theorem 12 of the  paper \cite{FrSab},  we  can  state  the following:
\begin{theorem} \label{th2.4}   Let $\l\in \C\setminus [0,\infty)$. The  ${\mathfrak S}_p$-norms of the operator
$W_1{\frak F}_0(\lambda)W_2$   satisfy the estimates
$$
||W_1{\frak F}_0(\lambda)W_2||_{{\mathfrak S}_p}\leq C_q |\l|^{-1+3/(2q)}||W_1||_{L^{2q}}||W_1||_{L^{2q}},
$$
with $3/2\leq  q\leq 2$ and $p=2q/(3-q).$ In particular,
\begin{equation}\label{FrSab}
||W_1{\frak F}_0(\lambda)W_2||_{{\mathfrak S}_4}\leq C |\l|^{-1/4}\Bigl(\int_{\R^3}|V|^2dx\Bigr)^{1/2},
\end{equation}
where the positive   constant $C$ is independent of $V$  and $\l$.
\end{theorem}

Finally, we  are going to obtain an estimate for  the  Hilbert-Schmidt   norm  of  the operator $W_1{\frak F}(\l)W_2$.
 Let $\mu$ be  the  function
$$
\mu(\lambda,x,y)={e^{-i{\pi\/4}}\/\sqrt{(4\pi)^3}}\int_0^\iy
e^{{i\/4t}|x-y|^2}\Bigl(e^{-i{t^3\/12}}-1\Bigr)e^{it\,\Lambda}{dt\/t^{3/2}},
$$
where $x,y\in {\mathbb R}^3$  and $\lambda\in \overline{\mathbb C}_+.$
\begin{proposition}
There  exists a  universal constant $C>0$  such that
 \[\lb{3.28}|\mu(\lambda,x,y)|\leq C{(1+|x|)^2(1+|y|)^2\/(1+|\l|)^{1/4}},\qquad \text{for\,\, all}\,\, x,y\in {\mathbb R}^3\, \text{and} \,\lambda\in \overline{\mathbb C}_+.\]
\end{proposition}

{\it Proof.}  First, note that  the  function $\mu(\lambda,x,y)$  is    bounded  by a  constant independent of  the variables $\l$, $x$ and $y$.
The  latter implies  that one needs to prove \eqref{3.28} only  for $|\l|>1.$  For this purpose we
  set $\beta=|\l|^{1/2}$ and
 write $\mu(\l,x,y)$  as  the sum  of  two integrals
\[\label{r=}
\begin{aligned}
\mu(\lambda,x,y)={e^{-i{\pi\/4}}\/\sqrt{(4\pi)^3}}\int_0^\beta
e^{{i\/4t}|x-y|^2}\Bigl(e^{-i{t^3\/12}}-1\Bigr)e^{it\,\Lambda}{dt\/t^{3/2}}+\\ {e^{-i{\pi\/4}}\/\sqrt{(4\pi)^3}}\int_\beta^\iy
e^{{i\/4t}|x-y|^2}\Bigl(e^{-i{t^3\/12}}-1\Bigr)e^{it\,\Lambda}{dt\/t^{3/2}}.
\end{aligned}
\]
The  second  integral   can be easily estimated  by $C/|\l|^{1/4}$, while the first  integral   equals
$$
{ie^{-i{\pi\/4}}\/\l\sqrt{(4\pi)^3}}\int_0^\beta
e^{it\l}\frac{d}{dt}\Bigl[{e^{{i\/4t}|x-y|^2}\Bigl(e^{-i{t^3\/12}}-1\Bigr)e^{-it(x_1+y_1)/2}\/t^{3/2}}\Bigr]dt+O(1/|\l|^{7/4}),
$$
as $|\l|\to\iy.$ Now we  observe that
$$
\Bigl| \frac{d}{dt}\Bigl[{e^{{i\/4t}|x-y|^2}\Bigl(e^{-i{t^3\/12}}-1\Bigr)e^{-it(x_1+y_1)/2}\/t^{3/2}}\Bigr] \Bigr|\leq C (t^{-1/2}|x-y|^2+|x_1|+|y_1|+1+t^{1/2})
$$
and integrate  the   expression in the right hand  side with respect to $t$   from  $0$ to $\b.$
As a result, we  will obtain that
$$
|\mu(\lambda,x,y)|\leq C (|\l|^{-3/4}|x-y|^2+(|x_1|+|y_1|+1)|\l|^{-1/2}+|\l|^{-1/4})\qquad  \text{for}\,\, |\l|>1,
$$
which definitely  implies that
$$
|\mu(\lambda,x,y)|\leq C |\l|^{-1/4}\Bigl(|x-y|^2+|x_1|+|y_1|+1\Bigr), \qquad  \text{for}\,\, |\l|>1.
$$
It remains to note that $|x-y|^2+|x_1|+|y_1|+1\leq  (1+|x|)^2(1+|y|)^2$.
The proof is completed. $\,\BBox$

\bigskip

\begin{corollary}\lb{cor3.2}
Let ${\mathfrak F}(\lambda)$  be the operator on $L^2({\mathbb R}^3)$  with the integral kernel $\mu(\lambda,x,y)$.
Let $W_1,\,W_2\in C^\infty_0({\mathbb R}^3)$ be two functions  such that $|W_1|=|W_2|$  and let $V=W_1W_2$.
Then
\[
||W_1{\mathfrak F}(\lambda)W_2||_{{\mathfrak S}_2}\leq \frac C{(1+|\l|)^{1/4}}\int_{{\mathbb R}^3}(1+|x|)^{4}|V(x)|\,dx,
\]
where $C$ is the same as in \eqref{3.28}
\end{corollary}

\begin{corollary} \lb{co4.3}There  exists  a universal    constant  $C>0$, such that
all  eigenvalues $\lambda\in {\mathbb C}\setminus \R$ of the operator $H$ are  situated in the  disk
$$
|\l|^{1/4}\leq C\Bigl(\int_{{\mathbb R}^3}(1+|x|)^4|V(x)|\,dx+\Bigl(\int_{\R^3}|V|^2dx\Bigr)^{1/2}\Bigr).
$$
Moreover, there  is a universal constant  $C_1>0$ such that  the condition
$$
\int_{{\mathbb R}^3}(1+|x|)^4|V(x)|\,dx+\Bigl(\int_{\R^3}|V|^2dx\Bigr)^{1/2}<C_1
$$ implies  that the  spectrum of $H$  coincides  with the real line ${\mathbb R}$.

\end{corollary}

{\it Proof.}  It is  very well known that all non-real  eigenvalues  of $H$ are  situated in  the set
$$
\{\l\in {\mathbb C}\setminus \R:\,\,\, ||Y_0(\lambda)||\geq 1) \}.
$$
Obviously,
 $$||Y_0(\lambda)||=||W_1R_0(\l)W_2||\leq||W_1{\frak F}_0(\l)W_2||_{{\mathfrak S}_4}+||W_1{\frak F}_1(\l)W_2||_{{\mathfrak S}_2}+||W_1{\frak F}(\l)W_2||_{{\mathfrak S}_2}.$$
On the other  hand,  according to  Theorem~\ref{th2.4}  combined with Corollary~\ref{corollary3.3},
$$
||W_1{\frak F}_0(\l)W_2||_{{\mathfrak S}_4}\leq C\frac1{1+|\lambda|^{1/4}}\Bigl[\Bigl(\int_{\R^3}|V|^2dx\Bigr)^{1/2}+\Bigl(\int_{\R^3}|V|^{3/2}dx\Bigr)^{2/3}\Bigr].
$$
Due  to Corollary~\ref{col3.5},  we also have
$$
||W_1{\frak F}_1(\l)W_2||_{{\mathfrak S}_2}\leq C\frac1{1+|\lambda|^{1/4}}\Bigl[\Bigl(\int_{\R^3}(1+|x_1|)^{3/2}|V|dx\Bigr)+\Bigl(\int_{\R^3}(1+|x_1|)^{3/4}|V|^{3/2}dx\Bigr)^{2/3}\Bigr].
$$
Finally, Corollary~\ref{cor3.2}  gives  us   the estimate
$$
||W_1{\frak F}(\lambda)W_2||_{{\mathfrak S}_2}\leq \frac C{(1+|\l|)^{1/4}}\int_{{\mathbb R}^3}(1+|x|)^4|V(x)|\,dx.
$$
Consequently,
$$
||Y_0(\lambda)||\leq \frac C{(1+|\l|)^{1/4}}\Bigl[\int_{{\mathbb R}^3}(1+|x|)^4|V(x)|\,dx+\Bigl(\int_{\R^3}|V|^2dx\Bigr)^{1/2}\Bigr].
$$
The latter  implies  both statements of Corollary~\ref{co4.3}. $\,\,\,\,\,\,\,\,\,\,\,\,\,\,\,\,\,\BBox$

\bigskip

As a consequence of the method, we obtain the following estimate with a very short expression  in the right hand side:

\begin{theorem}\lb{thm3.4}
Let ${\rm Im}\l\geq0$. Then  for any $p>11$,  there exists  a  positive constant  $C_p>0$ depending only on $p$  such that
\[
\lb{upperY0}
||Y_0(\l)||_{{\mathfrak S}_4}\leq\frac{ C_p}{1+|\l|^{1/4}}\Bigl(\int_{{\mathbb R}^3}(1+|x|)^p|V|^2 dx\Bigr)^{1/2}.
\]

\end{theorem}
{\it Proof.} It is  enough to note that
$$
\int_{{\mathbb R}^3}(1+|x|)^4|V(x)|\,dx+\Bigl(\int_{\R^3}|V|^2dx\Bigr)^{1/2}\leq C_p\Bigl(\int_{{\mathbb R}^3}(1+|x|)^p|V|^2 dx\Bigr)^{1/2}.
$$
$\,\,\,\,\,\,\,\,\,\,\,\,\,\,\,\,\BBox$

\bigskip

  Let  us prove that  that all eigenvalues of $H$ are  contained in a  disk of a finite radius  under  the condition $V\in L^{q/2}(\R^3)\cap L^\infty(\R^3)$  where $q<3$.
First, we   find $\ve$ so that $q=3/(1+\ve)$  and set $p=3/(1-\ve).$  Without  loss of generality,
we can assume that $\ve$  is very small.

According to \eqref{interpol5.40},
\[
||W_1R_0(\l)W_2||_{{\frak S}_p}\leq C\Bigl( ||W_1||_{L^p}||W_2||_{L^p}+||W_1||_{L^q}||W_2||_{L^q}\Bigr),\qquad \forall \l\in \C\setminus \R.\notag
\] Choose $\d>0$  so small  that $C\d  ( ||W_1||_{L^p}+||W_1||_{L^q}+||W_2||_{L^p}+||W_2||_{L^q}+2\d)<1/2$.
Let  now $\tilde W_j$ be bounded compactly supported approximations of $W_j$  having the property that
\[
 ||\tilde W_j-W_j||_{L^p}+||\tilde W_j-W_j||_{L^q}<\d.\notag
\]
Then
\[\lb{9.1}
||W_1R_0(\l)W_2-\tilde W_1R_0(\l)\tilde W_2||_{{\frak S}_p}\leq  C\d  \Bigl( ||W_1||_{L^2}+||W_1||_{L^3}+ ||W_2||_{L^2}+||W_2||_{L^3}+2\d\Bigr).
\]
The right hand  side in \eqref{9.1}   is  smaller than $1/2.$
Set  now $\tilde Y_0(\l)=\tilde W_1R_0(\l)\tilde W_2$. Then according to  the methods  that led us to  \eqref{upperY0},
\[
||\tilde Y_0(\l)||_{{\mathfrak S}_4}=O(|\l|^{-1/4}),\qquad \text{as}\quad |\l|\to\iy.
\]
It  remains to   note  that $||Y_0(\l)||\leq \frac12+||\tilde Y_0(\l)||$. $\BBox$

\bigskip

We provide an extensive  list of mathematical articles \cite{AF85}-\cite{CD78},  \cite{H77}-\cite{KP03},  \cite{{K86}}, \cite{LQZ89}- \cite{PushS}, \cite{RS76}, \cite{Y79},  \cite{Y81}  containing the   important  work on  Stark operators, which are  operators with the potential corresponding to  a  constant  electric  field,  and  the work related to the  study of the   Stark  effect.  Our  list  includes  the titles  of the books \cite{RS78}   and \cite{S05}  containing  the relevant theory  of
Schr\"odinger operators  and perturbation determinants.
Finally,  the  paper  \cite{EKtrace} is mentioned because it  indicates  the possible  direction  of the  related follow  up   research.

\end{document}